\documentclass[trsc,nonblindrev]{informs3}

\OneAndAHalfSpacedXI % current default line spacing
\usepackage{natbib}
 \bibpunct[, ]{(}{)}{,}{a}{}{,}%
 \def\bibfont{\small}%
\TheoremsNumberedThrough     % Preferred (Theorem 1, Lemma 1, Theorem 2)
\EquationsNumberedThrough    % Default: (1), (2), ...
\usepackage{ccaption}
\usepackage{amsmath}
\usepackage{subfigure}
\usepackage{empheq}
\makeatletter
\renewcommand{\fnum@figure}{\textbf{FIGURE~\thefigure} }
\renewcommand{\fnum@table}{\textbf{TABLE~\thetable} }
\makeatother
\captiontitlefont{\bfseries \boldmath}
\usepackage{cases}
\captiondelim{\;}
\usepackage{hyperref} % To include references inside the pdf file
\usepackage{amssymb}
\usepackage{graphicx}
\usepackage{rotating,tabularx}
\usepackage{tikz}
\newcommand{\R}{\mathbb{R}}
\usepackage{pstricks}
\usepackage{xstring}
\definecolor{darkblue}{rgb}{0,0,0.5}
\hypersetup{colorlinks=true,linkcolor=black,citecolor=darkblue,urlcolor=darkblue}
\newcommand{\newref}[2]{\hyperref[#2]{#1~\ref*{#2}}} % hyperref: \newref{labeltext (e.g. table, figure, ...)}{labelname}
\def\Bc{{\cal B}}
\def\Cc{{\cal C}}

\def\Fc{{\cal F}}

\def\Mc{{\cal M}}

\def\Pc{{\cal P}}

\def\Tc{{\cal T}}
\def\Fc{{\cal F}}
\def\Uc{{\cal U}}

\def\Xc{{\cal X}}

\def\1{{\bf 1}}

\def \N{\mathbb{N}}

\def\argmin{\mathop{\rm argmin}}
\def\argmin_#1{\underset{#1}{\mathrm{argmin\, }}}

\def \ep{\hbox{ }\hfill$\Box$}
\def \lip{\hbox{Lip}}

\begin{document}
%%%%%%%%%%%%%%%%

\setcounter{page}{1}
\VOLUME{00}%
\NO{0}%
\MONTH{Xxxxx}% (month or a similar seasonal id)
\YEAR{0000}% e.g., 2005
\FIRSTPAGE{000}%
\LASTPAGE{000}%
\SHORTYEAR{00}% shortened year (two-digit)
\ISSUE{0000} %
\LONGFIRSTPAGE{0001} %
\DOI{10.1287/xxxx.0000.0000}%

\RUNAUTHOR{Ameli et al.}

\RUNTITLE{Departure Time Choice Models Based on Mean Field Games}

% \TITLE{Bundling Information Goods of Decreasing Value}
\TITLE{Departure Time Choice Models in Urban Transportation Systems Based on Mean Field Games}

\ARTICLEAUTHORS{%
\AUTHOR{Mostafa Ameli}
\AFF{Univ. Gustave Eiffel, COSYS, GRETTIA\\
14-20 Boulevard Newton, 77420 Champs-sur-Marne, France, \href{mailto:mostafa.ameli@univ-eiffel.fr}{mostafa.ameli@univ-eiffel.fr}}%, \URL{}}

\AUTHOR{Mohamad Sadegh Shirani Faradonbeh}
\AFF{Department of Management sciences, University of Waterloo\\
200 University Avenue W, Waterloo, Ontario, Canada N2L 3G1, \href{mailto:msshiran@uwaterloo.ca}{msshiran@uwaterloo.ca}}%, \URL{}}

\AUTHOR{Jean-Patrick Lebacque}
\AFF{Univ. Gustave Eiffel, COSYS, GRETTIA\\
14-20 Boulevard Newton, 77420 Champs-sur-Marne, France, \href{mailto:jean-patrick.lebacque@univ-eiffel.fr}{jean-patrick.lebacque@univ-eiffel.fr}}%, \URL{}}

\AUTHOR{Hossein Abouee-Mehrizi}
\AFF{Department of Management sciences, University of Waterloo\\
200 University Avenue W, Waterloo, Ontario, Canada N2L 3G1, \href{mailto:haboueem@uwaterloo.ca}{haboueem@uwaterloo.ca}}%, \URL{}}

\AUTHOR{Ludovic Leclercq}
\AFF{Univ. Gustave Eiffel, Univ. Lyon, ENTPE, LICIT\\
3 Rue Maurice Audin, 69518 Vaulx-en-Velin cedex, France, \href{mailto:ludovic.leclercq@univ-eiffel.fr}{ludovic.leclercq@univ-eiffel.fr}}%, \URL{}}
} % end of the block

\ABSTRACT{%
Departure time choice models play a crucial role in determining the traffic load in transportation systems. This paper introduces a new framework to model and analyze the departure time user equilibrium (DTUE) problem based on the so-called Mean Field Games (MFGs) theory. The proposed framework is the combination of two main components including  (i) the reaction of travelers to the traffic congestion by choosing their departure times to optimize their travel cost; and (ii) the aggregation of the actions of the travelers, which determines the system level of service. The first component corresponds to a classic game theory model while the second one captures the travelers' interactions at the macroscopic level and describes the system dynamics. 
In this paper, we first present a continuous departure time choice model and investigate the equilibria of the system. Specifically, we demonstrate the existence of the equilibrium and characterize the DTUE. 
Then, a discrete approximation of the system is provided based on deterministic differential game models to numerically obtain the  equilibrium of the system. To examine the efficiency of the proposed model, we compare it with the departure time choice models in the literature. We apply our framework to a standard test case and observe that the solutions obtained based on our model are 5.6\% better in terms of relative cost compared to the solutions determined based on models in the literature. Moreover, our proposed model converges with less number of iterations than the reference solution method in the literature.  
Finally, the model is scaled-up to the real test case corresponding to the whole Lyon Metropolis with real demand pattern. The results show that the proposed framework is able to tackle much larger test case than usual to includes multiple preferred travel times and heterogeneous trip lengths more accurately than existing models in the literature.
}%

%\KEYWORDS{deterministic inventory theory; infinite linear programming duality; 
%  existence of optimal policies; semi-Markov decision process; cyclic schedule}

%\HISTORY{Submitted [\today]}

\maketitle
%%%%%%%%%%%%%%%%%%%%%%%%%%%%%%%%%%%%%%%%%%%%%%%%%%%%%%%%%%%%%%%%%%%%%%

% Samples of sectioning (and labeling) in TRSC
% NOTE: (1) \section and \subsection do NOT end with a period
%       (2) \subsubsection and lower need end punctuation
%       (3) capitalization is as shown (title style).
%
%\section{Introduction.}\label{intro} %%1.
%\subsection{Duality and the Classical EOQ Problem.}\label{class-EOQ} %% 1.1.
%\subsection{Outline.}\label{outline1} %% 1.2.
%\subsubsection{Cyclic Schedules for the General Deterministic SMDP.}
%  \label{cyclic-schedules} %% 1.2.1
%\section{Problem Description.}\label{problemdescription} %% 2.

% Text of your paper here

\section{Introduction}

In urban transportation systems, representing trip-making behavior requires a deep understanding of the interactions and interrelation between travelers' decisions and the systems' performance \citep{mahmassani1985dynamic}. In this context, studying the dynamics of travelers' departure time choice behavior, particularly in congested systems, is of fundamental importance \citep{ben2003discrete}. Departure time choice models represent how travelers choose their departure time considering their own desired arrival time \citep{hendrickson1981schedule}. The goal of each traveler in the system is to optimize its own travel cost \citep{dafermos1968traffic, sheffi1985urban}. This means that all travelers are assumed to be fully rational decision-makers. They anticipate other travelers' behavior in order to make an optimal decision, i.e., minimizing the additional travel cost due to the traffic congestion \citep{perakis2006analytical} and narrowing the final arrival time to the desired one \citep{smith1984stability, guo2018we}. Thus, the system behaves as a game in which the winners experience minimum travel cost \citep{cominetti2015dynamic} considering not only the travel cost but also the departure time. Based on the first principle of Wardrop \citep{wardrop1952road}, the game may have an equilibrium state called Departure Time User Equilibrium (DTUE) where no traveler can improve his individual travel cost by changing his departure time \citep{mahmassani1984dynamic, ran1993new}. 

Most of the previous studies on the DTUE problem in the literature are based on solving a classic Nash equilibrium problem coupled with %a classical traffic flow model at the link level \citep{wang2018dynamic, ban2012modeling}. For instance, the idea of DTUE has been widely used in the context of 
a single point-queue (bottleneck) model based on the pioneering paper of \cite{vickrey1969congestion}. The idea behind the point-queue models is to assume that the travel cost on the transportation system consists of a free-flow travel cost plus a congestion cost represented by a queueing cost \citep{daganzo1985uniqueness}. Therefore, DTUE arises because the queue capacity is limited, and travelers should consider the trade-off between the travel time losses and the costs corresponding to arriving later or earlier than the preferred arrival time \citep{ata2018equilibrium}. 
A comprehensive literature review on the bottleneck models has been recently conducted by \cite{li2020fifty}, which  highlights the developments and the applications of bottleneck models to transportation systems in the past half-century. 

Representing the urban transportation network by an origin, a destination, and a single bottleneck is not realistic \citep{lamotte2018morning, nagel2003still}. Indeed congestion depends on the detailed topology of the transportation network (e.g., the spatial distribution of origins, destinations, routes, and roads). In addition, congestion is impacted by the distribution of trips and vehicle densities \citep{jin2020generalized, ji2014empirical}. To integrate such features while keeping the macroscopic scale, a common is to represent the network dynamics by Macroscopic Fundamental Diagram (MFD) %approach in urban transportation economics and transportation science is to model the  transportation network using the bathtub 
models. This kind of model is also referred to bathtub models in the economic literature \citep{arnott2016equilibrium, mariotte2017macroscopic, lamotte2016morning}.  

\subsection{MFD/Bathtub models}
The first bathtub model was introduced by \cite{vickrey1991congestion, vickrey2020congestion}.  The travel demand is described by (i) the total number of trips, (ii) the distribution of trips' departure time, and (iii) the distribution of trips' length. MFD and Bathtub are two names for a set of equations corresponding to macroscopic traffic models. The single bathtub model considers an undifferentiated movement area to represent a dense network of congested links. The motion of travellers is assumed to take place at a speed which is considered to be uniform over the network but varies over time depending on the overall network loading \citep{bao2020leaving}. Therefore, the model does not need the location information of the origin and destination of travelers. %As in all macroscopic models, the speed is a function of the mean %transportation modes???%consider factoring out "transportation modes", replace by "mean"?%
%density in the area \citep{fosgerau2015congestion}.  
A trip is defined by its length and departure time in the dynamic setting. When a trip starts, its remaining distance to travel decreases following the evolution of the network mean speed. %The trip is started at the departure time when the trip length starts to decrease by the network supply (network speed). Afterward, the trip is finished when the remaining trip length becomes zero. 
Note that the network speed depends on the network characteristics (e.g., network size and road capacities) as well as the load on the network (network density) \citep{fosgerau2015congestion}. In order to capture the dynamics of the system, \cite{vickrey1991congestion} defined an ordinary differential equation to describe the evolution of the number of active trips (users in the network). Such a model resort to a strong assumption that the average remaining distance of active trips is constant. Another option is to assume that the remaining trip distance of active trips %constant ???, or the trip distance may be assumed is that it 
follows a time-independent negative exponential distribution \citep{vickrey1994, vickrey2019types}. However, based on the empirical studies of \cite{liu2012understanding, thomas2013empirical, tsekeris2013city}, travelers' trip length distribution is neither time-independent nor exponential. %Note that almost all macroscopic models such as classical bottleneck model \citep{vickrey1969congestion} and macroscopic fundamental diagram (MFD) \citep{daganzo2007urban} models have an assumption on the trip-length distribution. 

Recently, \cite{jin2020generalized} reviewed and analyzed several studies that relaxed the assumption on the trip length distribution \citep[see, e.g., ][]{leclercq2017dynamic, mariotte2017macroscopic, lamotte2018morning} and then proposed the generalized bathtub model which captures any distribution of the trip length. From a mathematical point of view, the total number of active trips is the  primary variable for most bathtub models. The generalized bathtub model focuses on the number of active trips with remaining distances greater than or equal to a threshold. By this definition, \cite{jin2020generalized} derived a set of partial differential equations to track the distribution of the remaining trip lengths. Further properties of the model are discussed in \cite{jin2020generalized}. In this work, we use the generalized bathtub model to capture the state of the urban transportation network for the departure time equilibrium problem.

%This study is focused on the calculation of the DTUE in this setting. \cite{arnott2018solving} developed a customized method for the computational solution of DTUE in the bathtub model.
%\cite{small2003hypercongestion} used the classical bathtub model to address the departure time user equilibrium with homogeneous travelers. They assumed that the travel time of each trip is determined by the conditions at its entrance. This assumption causes substantially different between experienced travel time and the instantaneous travel time.
\subsection{Departure time choice problem}

The departure time choice (also known as ``morning commute") problem at the network level is well-reviewed by \cite{lamotte2018congestion}. One of the main questions that have not been well studied in the literature is  how a departure time choice model can take into account the heterogeneity of the trip lengths with multiple preferred arrival time \citep{lamotte2018congestion}. The most complex equilibrium problem that has been addressed in the literature is modeling and numerically solving the DTUE problem for a group of travelers with a single distributed preferred arrival time and heterogeneous trip lengths \citep{lamotte2018morning}. This model is supported by empirical data and simulation. Note that even when a simpler bottleneck formulation is used for solving the DTUE problem (reviewed by \cite{jin2020stable}), few studies in the literature consider multiple preferred arrival time for commuters \citep{akamatsu2020new, lindsey2004existence, doan2011existence, ramadurai2010linear, takayama2017bottleneck, akamatsu2018departure, lindsey2019equilibrium}. Recall that the heterogeneity of the travelers' trip distance is not considered by single bottleneck models because they consider a single origin-destination \citep{akamatsu2020new}. The goal of this paper is to develop a more general mathematical framework to address departure time choice equilibrium with heterogeneous trip lengths and many desired arrival times in an urban transportation network.  

The concept of DTUE, originally, comes from game theory and Nash-equilibrium principles \citep{sun2017decision}. In general, with rational travelers, the user equilibrium problems represent fixed points \citep{wang2018dynamic, bortolomiol2019disaggregate}. The equilibration process of DTUE models is always addressed by population game theory at the network level \citep{arnott2013bathtub, yang2005evolutionary, arnott2018solving}. Population games have one strong assumption, which is called Myopia. Myopia in our problem means that travelers only take into account the current utilities of each alternative when choosing the departure time, without predicting other users' reactions (i.e., the departure time adjustment) \citep{sandholm2015population}. User interactions create new system states, which have new perceived utilities as a result of the evolution. This evolution process pushes the system at each iteration or day-to-day process toward equilibria. \cite{iryo2019instability} proved that when an evolution dynamics plays the role of the replicator dynamics, no stable equilibrium solution can be determined in the DTUE problem even when the demand profile is homogeneous, i.e., all users have the same travel distance. This study aims to overcome this limitation by employing a mean-field approximation and deriving a macroscopic framework. Note that this approximation does not mean that the users have perfect knowledge about the other users and the system. Because if travelers used perfect information, there is no need for an iterative or day-to-day process. Therefore, we deploy Mean field games to represent a prediction model of users regarding the other users decisions and evolution of the network.

\subsection{Mean field games}

To propose a new perspective on the DTUE, we resort to the Mean Field Games (MFGs) framework. The mathematical foundations of this theory were introduced in the seminal papers of \cite{lasry2006jeux, lasry2007mean}. The theory and methodology of MFGs have rapidly developed in different engineering fields \citep{djehiche2016mean}. The theory of MFGs studies decision-making  problems with an infinite number of interacting players \citep{adlakha2013mean}. The MFGs theory restates the classical game theory model as a micro-macro model \citep{cardaliaguet2010notes}. It allows defining players at the microscopic level similar to classical game theory models while translating the effect of players' decision to macroscopic models \citep{caines2015mean}. Therefore, instead of solving a large set of highly coupled equations that represent the interactions among players on a microscopic level, the core idea of MFGs is to exploit the ``smoothing" effect of large numbers of interacting players. The MFGs' main assumption (called mean field approximation) states that each player only reacts to a ``mass", which is defined by aggregating the effect of all the players. This approach simplifies the complex multi-agent dynamic systems at a macroscopic level \citep{degond2014large}. %An important class of models fitting into the MFGs framework is \textit{Congestion game}\citep{lacker2018mean}.   

There are a few studies in the literature that apply MFGs to analyze  transportation systems and most of them apply MFGs theory in the context of control theory \citep{chevalier2015micro, huang2019game}, vehicle routing problem \citep{tanaka2020linearly, salhab2018mean} or pedestrian moving models \citep{aurell2019modeling}. This paper, for the first time, develops a MFGs-based framework for the departure time equilibrium problem. %Notably, a departure time choice model is considered involving a large number of travelers. Thus, the mean field approximation can be applied to the problem. 
In our framework, each traveler looks for the optimal departure time by predicting the other travelers' departure time choices, given the current information of the traffic network congestion (mean-field), which is extracted from the generalized bathtub model. Then, the mean field is updated based on the optimal departure time choice of the travelers. The Nash equilibrium state occurs when the initial mean field approximation of the system is equal to the final mean field derived from the travelers' optimal departure time distribution. This process is equivalent to solving a fixed-point problem \citep{friesz1993variational}. 

%Regarding the calculation process of the DTUE, many studies limit the solution space of the DTUE points by applying strong assumptions on the trip-length distribution of the demand profile.

To numerically solve the DTUE model and determine the equilibrium of the system, many studies in the literature limit the feasible space of the DTUE problem by making strong assumptions on the trip-length distribution of the demand profile.
For instance, recent studies on the morning commute problem assume that the optimal solution fulfills some sorting property relating to the trip length and departure time, e.g., First-In, First-Out (FIFO) by \citep{daganzo2015distance}, partial FIFO by \cite{lamotte2018morning}, and Last-In, First-Out (LIFO) by \cite{fosgerau2015congestion}. Such assumptions restrict the exploration of the solution space. Moreover, most approaches %to deal with the rescheduling process (either optimization or day-to-day convergence) 
in the literature have a common drawback: they do not guarantee the solution's optimality while they are costly computationally at large-scale \citep{huang2020dynamic}. In this paper, we relax all the assumptions concerning sorting properties in the solution method (i.e., departure time rescheduling process) to better explore the solution space. We also propose a new heuristic method to speed-up the calculation process while converging to a solution which is closer  to the DTUE equilibrium compare to the existing methods in the literature.  

%Moreover, most solutions to deal with the rescheduling process (either optimization or day-to-day convergence) in the literature suffer from classical limitations of classical search algorithms and heuristic methods.

%Our main results are that the proposed MFGs framework can find fast a closer solution to the DTUE that does not have any regular sorting property by relaxing the mentioned assumptions on modeling and optimization. 
In this study, we first express the dynamic departure time choice problem at the network level based on the mean field games theory and generalized bathtub model (\newref{Section}{sec:MFGs}). Then, we discuss the properties of the DTUE in the continuous and discrete settings and prove that the model can  represent the morning commute problem without any strong assumptions of homogeneity on the demand profile (\newref{Sections}{Sec:3} and\newref{}{Sec:4}). Finally, in \newref{Section}{Sec:6}, we evaluate the performance of the model against one of the recently proposed approaches in the literature to solve the DTUE (\newref{Section}{Sec:61}) and apply the proposed model to the real test case of Lyon Metropolis network (\newref{Section}{Sec:62}). We numerically demonstrate that the model can not only consider heterogeneous demand profile for the morning commute problem but also a large transportation system with a high number of travellers. In  \newref{Section}{Sec:Con}, we provide concluding remarks.  

%%%%%%%%%%%%%%%%%%%%

%The rest of the paper is organized as follows. \newref{Section}{sec:MFGs} presents the DTUE mathematical model. The continuous Mean Field Games framework and the existence of the equilibrium solution are discussed in \newref{Section}{Sec:3}. The complete MFGs model in a discretized setting and the solution method to solve the DTUE problem are presented in \newref{Section}{Sec:4}. \newref{Section}{Sec:6} presents the test cases for validating the model and the results obtained from the experiments on the large-scale network. \newref{Section}{Sec:Con} states concluding remarks and introduces future directions of work.  

%%%%%%%%%%%%%%%%%%%%%%%%%%%%%%%%%%%%%%%%%%%%%%%%%%%%%%%%%%%%%%%%%%%%%%%%%%%%%%%%%%%%%%%%%% 

\section{Problem Definition}\label{sec:MFGs}

\begin{table}[h] 
\begin{center}
\caption{List of notations}
\begin{tabular}{p{4cm}p{12.5cm}} \hline
% \multicolumn{2}{c}{Summary of Notations.} \\
% \textit{Basic notations:} \\
% $n$                 &   A positive integer number, $n \in \N$.\\
% $[n]$               &   Set of all positive integers less than or equal to $n$.\\
% $x^i$               &   The $i$-th component of vector $x, \ \forall x \in \R^n$, $i=1,\ldots, n$. \\
% $x^{-i}$            &   Sub-vector of vector $x$ excluding its $i$-th component. \\
% $C$                 &   A compact set of $\R^n$.\\
% $\Cc(C)$            &   The space of all real valued continuous functions on $M$.\\
% $\Pc(C)$            &   The space of all Borel probability measures on $M$.\\
% $d(u, v)$           &   The uniform distance between $u$ and $v$ on $\Cc(M)$.\\
% % $W_1(F,G)$          &   Wasserstein distance between $F$ and $G$ on $\Pc(M)$.\\
% \vspace{0.1cm}
% \textit{MFGs model:} \\
% $N$                 &   Set of independent trips.\\
$\Tc$                 &   Time horizon.\\
$n$                 &   The total number of trips.\\
$i$                 &   Index of trips, $i \in N$.\\
% $x$                 &   Vector of trip lengths. %(given by the demand profile).
$X_{max}$           &   Maximum trip length.\\
$X_{min}$           &   Minimum trip length.\\
% $t_a$               &   Vector of desired arrival times. %(given by the demand profile).
$x^i$               &   Trip length of trip $i$.\\
$t^i_d$             &   Departure time of trip $i$.\\
% $t^{-i}_d$          &   Distribution of departure times of all trips excluding trip $i$.\\
$T(t^i_d, x^i)$     &   Travel time of a trip started at $t^i_d$ with trip length $x^i$.\\
$t^i_a$             &   Desired arrival time of trip $i$.\\
$\bar{t}^i_a$       &   Actual arrival time of trip $i$.\\
$v_t$               &   Velocity of the system at time $t$.\\
$c_t$               &   Fraction of the total demand that traveling in the system at time $t$.\\
$z(t)$              &   Characteristic travel distance.\\
$o_t$               &   Outflow fraction of the system at time $t$.\\
$\varphi(t,\cdot)$      &   Probability density function of the active trips’ remaining distances at $t$.\\
$\Phi(t,x)$         &   Fractions of active trips with trip lengths more than $x$ at time $t$.\\
$F$            &   In-flow measure,  the empirical distribution of the departures.\\
$V$        &   Speed function, which maps the fraction of active travellers to the velocity.\\

$\Delta t$ & Small time interval. \\
$\Delta x$ & Small space interval. \\

\hline \label{TAb:1}
\end{tabular}
\end{center} 
\end{table}

\newref{Table}{TAb:1} presents the list of notations used in this paper. Consider a system with $n$ independent trips indexed by $i \in [n]:=\{1,2,\ldots,n\}$ in a time horizon $\Tc := [0,T_{max}]$. The trip length of the $i$-th trip is denoted by $x^i \in \Xc:= [X_{min}, X_{max}]$. The goal of the player $i$ is to choose his departure time $t_d^i \in \Tc_d$ to arrive at the desired arrival time $t_a^i \in \Tc_a$, where $\Tc_d$ and $\Tc_a$ are two compact subsets of $\Tc$. Assume that the joint distribution of the desired arrival times and trip lengths given as the demand profile $m$. If we define
\begin{align}
    \label{eq:demand_profile}
    m_n := \frac{1}{n} \sum_{i=1}^n \delta_{(x^i,t_a^i)},
\end{align}
where $\delta$ denotes the Dirac delta function, then $m_n \Rightarrow m$ as $n \rightarrow \infty$. Here $ \Rightarrow$ shows the weak convergence of the measures. That means, $\int_{\Tc_d\times \Xc} \phi dm_n \rightarrow \int_{\Tc_d\times \Xc} \phi dm$ for all $\phi\in \Cc_b(\Tc_d\times \Xc)$, the set of all bounded continuous functions on ${\Tc_d\times \Xc}$ (see \cite{billingsley2013convergence,carmona2018probabilistic}). Note that $m$ is a \textit{probability measure}\footnote{Any function with values in $[0,1]$, returning $1$ for the entire space and $0$ for the empty set that satisfies countable additivity property. A function $F$ is countable additive if for all countable family $\{B_i\}$ of pairwise disjoint sets, it holds true that $F(\bigcup B_i)=\sum F(B_i)$, see \cite{billingsley2012probability}.} in the space of all probability measures defined on $\Xc \times \Tc_a$, i.e. $m \in \Pc(\Xc \times \Tc_a)$. Therefore, it fully describes the demand characteristics. Hereafter, we assume that $\{t_a^i,x^i\}_{i=1}^n$ are i.i.d random variables with the distribution $m$. We also make the following regularity assumption on $m$.

\begin{assumption}
\label{asm:regularity}
There exists a constant $M_m$ such that
\begin{equation}
    \label{eq:m-regulatity}
    m(B,\Tc_a) \leq M_m \lambda (B), \quad \forall B\in \Bc (\Xc)
\end{equation}
where $\Bc(\Xc)$ and  $\lambda$ denote the $\sigma$-algebra of Borel sets and the Lebesgue measure on $\Xc$, respectively.
\end{assumption}
Broadly speaking, Assumption \ref{asm:regularity} means that the demand is spread smoothly over $\Xc$ and particularly that the demand cannot be concentrated in Dirac distributions.

The congestion in the system at time $t$ is defined by the fraction of the total demand that is active at time $t$, which is captured by,
\begin{equation}
    \label{eq:congestion}
    c_t = \frac{1}{n} \sum_{i=1}^{n} \1_{[t_d^i, \bar{t}_a^i)} (t),
\end{equation}
where $\bar{t}^i_a$ denotes the actual arrival time of the $i$-th player and $\1_{[t_d^i, \bar{t}_a^i)}$ is an indicator function which returns $1$ if $t \in [t_d^i, \bar{t}_a^i)$ and $0$ otherwise. 
We assume that the velocity of the system at time $t$ depends on the fraction of travelling users in the system $c_t$ which is defined by a strictly decreasing speed function $V: \R^+ \mapsto \R^+$. Therefore, $V$ represents the mean network speed and is the key collective behavioral characteristic of the generalized bathtub model \citep{jin2020generalized}. Recall that the generalized bathtub model considers the network characteristics and the network load to calculate the network speed. We simply denote the velocity at time $t$ by $v_t := V(c_t)$ and assume that the velocity is the same for all players who are travelling at the same time.

To determine the travel time of a player, we first define a virtual user who starts his trip at time $0$. Then, the characteristic travel distance $z(t)$, travelled by this virtual user up to time $t$, is
\begin{align}
    \label{eq: virtual user}
    z(t) := \int_0^t v_s ds = \int_0^t V\left(c_s\right) ds.
\end{align}
%Also, define $z(t)$ as the distance travelled by a virtual user up to time $t$. Then, we have:
%\begin{align}
 %   \label{eq: virtual user}
  %  z(t) := \int_0^t v_s ds,
%\end{align}
Since $v_t >0\ \forall t\in \Tc$,
$z$ is an invertible function. Let $z^{-1}$ denote the inverse function of $z$. Then, we have $z^{-1}\big(z(t)\big) = t$ and $z^{-1}(x)$ represents the time at which the virtual user has reached $x$.

Now, let $T(t_d^i,x^i)$ denote the travel time of a player departing at time $t_d^i$ with trip length $x^i$. Considering (\ref{eq: virtual user}),  $T(t_d^i,x^i)$ can be determined by,
\begin{align}
    \label{eq: trip time}
    T(t_d^i,x^i) = z^{-1}\big(x^i+z(t_d^i)\big) - t_d^i.
\end{align}

%Accordingly, 
%optimizing the individual departure time choice, given the departure times and lengths of the other trips,  player $i$ would be able to calculate his travel time $T(t_d^i,x^i)$. Subsequently, using Equation (\ref{eq:congestion}),  (\ref{eq: trip time}), and the congestion-velocity relation $v_t = V(c_t)$, the actual arrival time $\bar{t}_a^i$ can be obtained. 
To determine the optimal departure time, we assume that each player aims to minimize his travel cost. In the DTUE problem, the travel cost is usually defined based on $\alpha$-$\beta$-$\gamma$ scheduling preferences \citep{fosgerau2015congestion}. That means, the cost function is defined as the sum of the travel time and a penalty cost for arriving at $t_d^i + T(t_d^i,x^i)$ instead of the desired arrival time. Specifically, we assume that each player's cost function is given by, 
\begin{align}
\label{eq:cost function}
J_i(t_d^i,t^i_a; t_d^{-i},x^{-i}) = \alpha T(t_d^i,x^i) + \beta\big(t_a^i - t_d^i - T(t_d^i,x^i)\big)_+ + \gamma\big(t_d^i + T(t_d^i,x^i) - t_a^i\big)_+,
 % Mod JPL   J_i(t_d^i; t_d^{-i}) = \alpha T(t_d^i,x^i) + \beta\big(t_a^i - t_d^i - T(t_d^i,x^i)\big)_+ + \gamma\big(t_d^i + T(t_d^i,x^i) - t_a^i\big)_+,
    % \\
    % J_i(t_d^i; t_d^{-i}) = T(t_d^i,x^i) + \frac{\alpha}{2}\big(t_d^i + T(t_d^i,x^i) - t_a^i\big)^2.
\end{align}
where $\alpha$ denotes the cost of travelling per unit of time, $\beta$ and $\gamma$ denote, respectively, the cost of earliness and lateness for the traveler arrival. Note that $(y)_+ = \max\{y,0\}$ as well $t_d^{-i}$ and $x^{-i}$ respectively express the dependency of $J$ on the departure times and trip lengths of the other users ($\not= i$) via their travel times.

The cost function defined in (\ref{eq:cost function}) captures the fact that travelers prefer not to deviate from their desired arrival time (i.e., arrive as close as possible to their desired arrival time) while they do not spend too much time on the traffic.
Note that the dependency of the cost function on
% the desired arrival times and 
the trip lengths is not emphasized in the notation, while it holds implicitly. 
%Provided that the velocity  is known throughout the time horizon, by Equation (\ref{eq: trip time}), the travel time can be determined and thus $t_d^i$ and $\bar{t}_a^i$ have a one-to-one relation. Therefore, given either one, the other one can be calculated.
Below, we provide the definition of the optimal strategy that each player adopts to determine his departure time.
\begin{definition}
\label{def:NE}
The departure time vector $\hat{t}_d := (\hat{t}_d^1,\ldots,\hat{t}_d^n) \in \Tc_d^n$ is a Nash equilibrium (NE) for the cost function given in (\ref{eq:cost function}), if for all $i \in [n]$ we have
\begin{align}
    \label{eq:nash-eq}
    J_i(\hat{t}_d^i,t_a^i; \hat{t}_d^{-i},x^{-i}) \leq J_i(t,t_a^i; \hat{t}_d^{-i},x^{-i}), \quad \forall t \in \Tc_d.
%  Mod JPL     J_i(\hat{t}_d^i; \hat{t}_d^{-i})
%    \leq 
 %   J_i(t; \hat{t}_d^{-i}), \quad \forall t \in \Tc_d.
\end{align}
\end{definition}
The above definition indicates that at a NE point $\hat{t}_d$, no player can decrease his travel cost by deviating from his departure time.
Based on Definition \ref{def:NE}, we define DTUE as a NE of the following Departure Time Choice Problem (DTCP):
% \begin{equation}
% \label{eq: DTCP}
%     \tag{DTCP}
 %  Mod JPL       \min_{t_d^i \in \Tc_d} J_i(t_d^i; t_d^{-i}) &= \alpha T(t_d^i,x^i) + \beta\big(t_a^i - t_d^i - T(t_d^i,x^i)\big)_+ + \gamma\big(t_d^i + T(t_d^i,x^i) - t_a^i\big)_+
%  \label{eq: DTCP}
 \begin{align}
    % \label{eq: DTCP}
 \tag{DTCP}
 \label{eq: DTCP}
 \min_{t_d^i \in \Tc_d} J_i(t_d^i,t_a^i; t_d^{-i}, x^{-i}) &= \alpha T(t_d^i,x^i) + \beta\big(t_a^i - t_d^i - T(t_d^i,x^i)\big)_+ + \gamma\big(t_d^i + T(t_d^i,x^i) - t_a^i\big)_+ \quad \forall i 
        \\
        \notag
        \textit{s.t.}\\
        &\begin{cases}
        % \tag{DTCP-Constraints}
        \label{eq:DTCP-C}
        c_t = \frac{1}{n} \sum_{j=1}^{n} \1_{\big[t_d^j, t_d^j + T(t_d^j,x^j)\big)} (t),
        \\
        z(t) = \int_0^t V(c_s) ds,
        \\
        T(t_d^i,x^i) = z^{-1}\big(x^i+z(t_d^i)\big) - t_d^i.
        \end{cases}
\end{align}
% \end{equation}
Similar to (\ref{eq:cost function}), \ref{eq: DTCP} provides the cost function of the $i$-th player. Note that, given the departure times and trip lengths of others ($t_d^{-i},\ x^{-i}$), the player $i$ is able to find his travel time. Specifically, according to the set of equations given in (\ref{eq:DTCP-C}), one can derive the characteristic travel distance $z(t)$ based on the fraction of active trips $c_t$. Then, the travel time function $T$ can be obtained.

Since analyzing the \ref{eq: DTCP} for a large $n$ is arduous due to ``curse of dimensionality", in the next section, we examine the behaviour of players in a system where the number of players goes to infinity, i.e., $n \rightarrow \infty$. This means that we adopt the MFGs approach to determine the DTUE.

\section{Mean Field Games Framework} \label{Sec:3}

In this section, we discuss the \ref{eq: DTCP} in the framework of the MFGs. First, recall that the idea behind the MFGs is to consider a proxy that represents the macroscopic behavior of all the players at once, instead of taking into account their departure times individually. Therefore, to capture the information of entering trips from the viewpoint of the $i$-th player when there are $n$ players (including player $i$) in the game, we define the following empirical measures,
\begin{align}
    \label{eq: in-flow1}
% Insert JPL    F^{i}_n =  \frac{1}{n-1} \sum_{j\in[n]\setminus i} \delta_{t_d^j, x^j}, 
    %  \begin{array}{l}
         &F^{i}_n :=  \frac{1}{n-1} \sum_{j\in[n]\setminus i} \delta_{t_d^j, x^j}, \\
         \label{eq: dis-in-flow1}
         &E^{i}_n :=  \frac{1}{n-1} \sum_{j\in[n]\setminus i} \delta_{t_d^j, x^j,t_a^j},
    %  \end{array}
\end{align}
where $\delta$ is the Dirac delta function. Note that the cost function of each player, defined in (\ref{eq:cost function}), is a symmetric function\footnote{For any $i \in [n]$, the cost function of the player $i$ satisfies,
\begin{align*}
    J_i(t_d^i, t_a^i; t_d^{-i}) = J_{\zeta(i)}(t_d^{\zeta(i)}, t_a^{\zeta(i)}; t_d^{-\zeta(i)}),
\end{align*}
for all permutation $\zeta$ on $\{1,\ldots,n\}$, see p. 49 of \cite{lacker2018mean}}. That means, the objective function of a player does not change if other players change their labels with each other. Further, as $n \rightarrow \infty$ the impact of each player on the system vanishes. That means by changing either the departure time, desired arrival time, or trip length of a player, the velocity of the system would be left unchanged. Thus, we define the in-flow measure $F$ and the dis-aggregated in-flow measure $E$ as the limits of the sequences $\{F^{i}_n\}_{n \in \N}$ and $\{E^{i}_n\}_{n \in \N}$, respectively. That is,
\begin{align}
    \label{eq: in-flow}
% Insert JPL       F = \lim_{n\rightarrow \infty} F^{i}_n.
    %   \begin{array}{l}
           &F := \lim_{n\rightarrow \infty} F^{i}_n,  \\
           \label{eq: dis-in-flow}
            &E := \lim_{n\rightarrow \infty} E^{i}_n.  
    %   \end{array}
\end{align}
Note that the limits are independent of $i$ as the impact of a single player vanishes as $n$ gets large. Further, $F$ is a probability measure on the product space $\Tc_d \times \Xc$, i.e., $F \in \Pc(\Tc_d \times \Xc)$, the set of all probability measures define on $\Tc_d \times \Xc$. In fact, for all $n \in \N$ and $i\in [n]$ the function $F^{i}_n$ is a probability measure and $F$ is the limit in the weak convergence sense.  % insert JPL 1
Similarly, one can show that
% $E$ is a probability measure on the product space $\Tc_d \times \Xc \times \Tc_a$, i.e., 
$E \in \Pc(\Tc_d \times \Xc \times \Tc_a)$. % JPL check \Tc_a 
According to (\ref{eq: in-flow1}) and (\ref{eq: dis-in-flow1}), the in-flow measure $F$ depends on $E$ in the following sense, 
% \begin{equation}
%     \label{eq:FdependentG}
%     \left< F,\varphi  \right> = \left< E,\phi \otimes  1_{\Tc_a} \right> \quad \forall \phi \in \Cc (\Tc_d \times \Xc ),
% \end{equation}
\begin{align}
    \label{eq:FdependentE}
    \int_{\Tc_d \times \Xc} \phi dF = \int_{\Tc_d \times \Xc \times \Tc_a} \phi \otimes 1_{\Tc_a} dE,
    \quad \forall \phi \in \Cc_b (\Tc_d \times \Xc ),
\end{align}
% where $\Cc (\Tc_d \times \Xc)$ denotes the set of all continuous functions on $\Tc_d \times \Xc$ and
where $\phi \otimes 1_{\Tc_a}$ is the tensor product of $\phi$ and $1_{\Tc_a}$. In fact,  (\ref{eq:FdependentE}) assures that the in-flow measure $F$ is the marginal probability measure of the departure times and trip lengths wrt to the dis-aggregated in-flow measure $E$, almost surely (a.s.). 
This linear dependency is continuous by (\ref{eq:FdependentE}) and is denoted by $\Fc: \Pc(\Tc_d \times \Xc \times \Tc_a) \rightarrow \Pc(\Tc_d \times \Xc)$ such that,
\begin{equation}\label{eq:FfunctionE}
    F = \Fc (E).
\end{equation}
Further, $E$ is constrained by the demand profile $m$, see (\ref{eq:demand_profile}), such that,
\begin{align}
    \label{eq:Edependentm}
    \int_{\Tc_d \times \Xc \times \Tc_a} 1_{\Tc_d}\otimes\phi dE = \int_{\Xc \times \Tc_a} \phi dm,
    \quad \forall \phi \in \Cc_b (\Xc \times \Tc_a).
\end{align}
Constraint (\ref{eq:Edependentm}) restricts the dis-aggregated  in-flow measure $E$ to a subset of $\Pc(\Tc_d \times \Xc \times \Tc_a)$ with marginal probability measures which is a.s. equal to the demand profile $m$.
% \begin{equation}\label{eq:DefConstrintsG}
%     \left< m,\varphi  \right> = \left< E, 1_{\Tc_d}  \otimes \phi \right> \quad \forall \phi \in \Cc ( \Xc \times \Tc_a ).
% \end{equation}
% The set of constraints applying to $G$ is denoted by $\Kc \subset \Pc(\Tc_d \times \Xc \times \Tc_a) $ and is defined by,
% \begin{equation}\label{eq:DomainK}
%     E \in \Kc \Longleftrightarrow \;\; \left[ \;
%     \begin{array}{l}
%       E \geq 0   \\
%       E  \;\; {\rm satisfies}\; {\rm equation  } \; (\ref{eq:DefConstrintsG}) 
%     \end{array}  \right.  
% \end{equation}
% End insert JPL 1
Constraints (\ref{eq:FdependentE}) and (\ref{eq:Edependentm}) together yield to the following demand constraint,
\begin{align}
    \label{eq: demand profile = marginal dist-x}
    F(\Tc_d,B) = m(B, \Tc_a), \quad \text{(a.s.)}
\end{align}
for all Borel measurable subsets of $\Xc$ such as $B$. Roughly speaking, from the demand viewpoint, the fraction of players having a trip length in $B$, $m(B, \Tc_a)$, matches the fraction of the departures with a trip length in $B$, $F(\Tc_d,B)$. Thus all the demand is served over the overall time period.

% To state and analyze the \ref{eq: DTCP} in MFGs framework, we need first to determine, similar to (\ref{eq:DTCP-C}), the relation between the travel time function $T$ and the in-flow measure $F$, i.e., the relation between the travel time of a player and others' departures. 

Now, given the in-flow measure $F$, an arbitrary player with the departure time $t_d$, trip length $x$ and desired arrival time $t_a$ can reconsider his cost criteria as a function $J: \Tc_d \times \Xc \times \Tc_a \times \Pc(\Tc_d \times \Xc) \mapsto \R^+$ defined as follows,
\begin{align}
    \label{eq: MFGs-cost function}
    J(t_d; x, t_a; F) = \alpha T(t_d,x) + \beta\big(t_a - t_d - T(t_d,x)\big)_+ + \gamma\big(t_d + T(t_d,x) - t_a\big)_+.
\end{align}
To be more specific, assuming $F$ is known, each player tries to choose his departure time $t_d$ by minimizing (\ref{eq: MFGs-cost function})\footnote{Obviously, each player knows his trip length $x$ and desired arrival time $t_a$.}. Towards this end, we should first state the travel time of a player with departure time $t_d$ and trip length $x$, i.e., $T(t_d,x)$, in terms of the in-flow measure $F$. That means, we want to find the relation between the travel time of a player and others' departures and trip lengths.

\subsection{Dynamics of the System}
\label{sec:SystemDynamics}
We assume that the in-flow measure $F$ with the probability density function $f$ is given. The goal of this section is to derive the dynamics of the characteristic travel distance $z$, defined in (\ref{eq: virtual user}). Then, using (\ref{eq: trip time}), we can clarify the relation between the travel time function $T$ and the characteristic travel distance $z$ which completes the definition of (\ref{eq: MFGs-cost function}).
% \begin{Assumption}
%     \label{asm: lipschitz of speed}
%     There is a Lipschitz continuous function $V: \R^+ \mapsto \R^+$ maps the fraction of travelling users $c_t$ to the traffic speed $v_t$, i.e. $v_t = V(c_t)$. That means there exists a constant $U>0$ such that:
%     \begin{align}
%         |V(x) - V(y)| \leq U|x-y|,\quad \forall x,y \in \R.
%     \end{align}
% \end{Assumption}

Assuming the number of players in the system goes to infinity ($n \rightarrow \infty$), we define the dynamics of the system according to the fraction of active trips instead of the number of them. Therefore, we denote by $\varphi(t,\cdot)$ the probability density function of the remaining trip lengths of active trips at time $t$. Then, the out flow of the system at time $t$, $o_t$, can be stated as follows,
\begin{align}
    \label{eq: outflow-explicit}
    o_t =  c_t v_t \varphi(t,0),
\end{align}
where, in the sense of (\ref{eq:congestion}), the system congestion in the limit model can be defined as
\begin{align}
    \label{eq:mfg_congestion}
    c_t = \lim_{n \rightarrow \infty} \frac{1}{n} \sum_{i=1}^{n} \1_{[t_d^i, \bar{t}_a^i)} (t)
\end{align}
To derive the dynamics of $\varphi(t,x)$, we use the idea of the generalized bathtub model, see \cite{jin2020generalized}. Note that in a system with $n$ trips, for a small time interval $\Delta t$, the number of active trips at $t+\Delta t$ with a remaining trip length in $[x, x+\Delta x]$ is $nc_{t+\Delta t}\varphi(t+\Delta t, x)\Delta x$. On the other hand, it is equal to the sum of new departures and trips with remaining trip length in $[x+v_t\Delta t, x+v_t\Delta t+\Delta x]$ at time $t$. Thus,
\begin{align*}
    &nc_{t+\Delta t}\varphi(t+\Delta t, x)\Delta x
    \approx
    nf(t,x)\Delta x\Delta t + nc_t\varphi(t,x+v_t\Delta t)\Delta x,
\end{align*}
which is equivalent to,
\begin{align*}
   &c_{t+\Delta t}\varphi(t+\Delta t, x)
    \approx
    f(t,x)\Delta t + c_t\varphi(t,x+v_t\Delta t).
\end{align*}
To simplify the system dynamics, we  approximate $c_{t+\Delta t}$ and $\varphi(t,x+v_t\Delta t)$ with $c_t + c'_t\Delta t$ and $\varphi(t,x) + v_t \frac{\partial_x \varphi(t,x)}\Delta t$, respectively. Then, dividing both sides by $\Delta t$ and letting   $\Delta t$ goes to zero, we get,
\begin{align}
    \label{eq: dynamic of varphi}
    c_{t}\partial_t\varphi(t,x) + c'_{t}\varphi(t, x)
    - c_tv_t\partial_x\varphi(t,x) = f(t,x).
\end{align}
We use $\partial_t$ and $\partial_x$ to denote, respectively, the partial derivative with respect to time $t$ and space $x$. Integrating both sides of (\ref{eq: dynamic of varphi}) with respect to $x$ from $x$ to $X_{max}$, we get,
\begin{align}
    \label{eq: dynamic of varphi-integrated}
    c_{t}\partial_t\Phi(t,x) + c'_{t}\Phi(t, x)
    - c_tv_t\partial_x\Phi(t,x) = \int_x^{X_{max}}f(t,\xi)d\xi.
\end{align}
Note that $\Phi(t,x)$ denotes the fraction of active trips with the remaining trip lengths more than $x$ at time $t$,
$$
   \Phi(t,x) = \int_x^{X_{max}} \varphi(t,\xi) \, d\xi \;. 
$$
Then, using relation $z'(t) = v_t$ as a direct result of (\ref{eq: virtual user}), we have,
\begin{align*}
    \frac{d}{dt} \big(c_{t}\Phi(t,x - z(t)) \big) = c_{t}\partial_t\Phi(t,x - z(t)) + c'_{t}\Phi(t, x - z(t))
    - c_tv_t\partial_x\Phi(t,x - z(t)).
\end{align*}
Therefore, the equality given in (\ref{eq: dynamic of varphi-integrated}) can be written as,
\begin{align*}
    \frac{d}{dt} \big(c_{t}\Phi(t,x - z(t)) \big) = \int_{x - z(t)}^{X_{max}}f(t,\xi)d\xi.
\end{align*}
Thus, by integrating both sides with respect to time from $0$ to $t$ and setting $y = x - z(t)$, the dynamics of the system can be presented as,
\begin{align}
    \label{eq: dynamics of the Phi}
    c_{t}\Phi(t,y)= \int_0^t \int_{y+z(t)-z(s)}^{X_{max}} f(s,\xi)d\xi ds.
\end{align}
Note that we assume the system is empty at time $0$. Moreover, taking partial derivative with respect to $x$ from both sides, considering that $c_t$ is independent of $x$, and applying Leibniz's integral rule, we obtain,
\begin{align}
    \label{eq:eq1}
    \partial_x \big(c_{t}\Phi(t,x)\big)
    =
    c_{t} \partial_x\Phi(t,x)
    =
    - \int_0^t f\big(s,x+z(t)-z(s)\big)ds.
\end{align}
Finally, using (\ref{eq: outflow-explicit}) and the definition of in-flow measure $F$ and its probability density function $f$, the dynamics of the fraction of the active trips at time $t$, $c_t$ satisfies
\begin{align*}
    c'_t = \int_0^{X_{max}} f(t,x)dx
    - o_t 
    = \int_0^{X_{max}} f(t,x)dx + c_t v_t \partial_x\Phi(t,x)|_{x=0}.
\end{align*}
Here, $\partial_x\Phi(t,x)|_{x=0}$ demonstrates the right derivative at $0$ as the left derivative is not defined. Substituting (\ref{eq:eq1}), gives,
\begin{align*}
    c'_t = \int_0^{X_{max}} f(t,x)dx - v_t \int_0^t f\big(s,z(t)-z(s)\big)ds.
\end{align*}
Integrating both sides with respect to time and using
Tonelli's theorem (see, e.g., Theorem 18.3 of \cite{billingsley2012probability}), to change the order of the integration, we get,
\begin{equation}
    \label{eq: dynamics of the system}
    \begin{aligned}
    c_t &= \int_0^t \int_0^{X_{max}} f(r,x)dxdr - \int_0^t  \int_0^r v_r f\big(s,z(r)-z(s)\big)ds dr
    \\
    &= \int_0^t \int_0^{X_{max}} f(r,x)dxdr - \int_0^t  \int_s^t v_r f\big(s,z(r)-z(s)\big)dr ds.
    % \\
    % &= \int_0^t F(dr, \Xc) - \int_0^t F\Big(ds,\big[0,z(t)-z(ds) \big]\Big)
    % \\
    % &= \int_0^t F\Big(ds,\big(z(t)-z(ds),X_{max} \big]\Big).
\end{aligned}
\end{equation}

In the light of the equality given in (\ref{eq: virtual user}) and considering that $z'(t) = v_t$, the result of the above discussion about the characteristic travel distance is summarized in Proposition \ref{prp: congestion dynamics}. To state the proposition rigorously, we make the following assumption. 
\begin{assumption}
    \label{asm: lipschitz of F}
    Let $G>0$ be a constant. Then, for all Borel measurable subset of $\Tc_d \times \Xc$ such as $B$, we assume that the in-flow measure $F \in \Pc(\Tc_d \times \Xc)$
    % \footnote{ $B \in \Bc(\Tc_d \times \Xc)$ where $\Bc(\Tc_d \times \Xc)$ is a Borel $\sigma$-algebra defined on $\Tc_d \times \Xc$, see  \cite{billingsley2012probability}.}
    satisfies
    \begin{align}
        \label{eq: lipschitz of F}
        F(B) \leq G \lambda_2(B),
    \end{align}
    where $\lambda_2$ is Lebesgue measure on $\R^2$.
\end{assumption}
Assumption \ref{asm: lipschitz of F} is a regularity condition. It means that the in-flow measure $F$ is spread smoothly both with respect to departure time and trip length. Assumption \ref{asm: lipschitz of F} is a technical condition which essential for most proofs of this section. 
% It means we will not have a concentration of departure time or trip lengths concurrently, i.e., the distribution of departure time and trip lengths are smooth. This assumption is reasonable in the concept of DTUE because it demonstrates that all users are not starting simultaneously, which can also generate an extra cost of congestion for users objective function. Besides, this assumption distinguishes the morning commute problem and disruption or evacuation models. 

% ??? make the following a full fledged definition?
Given the demand profile $m$, defined in (\ref{eq:demand_profile}), and the constant $G>0$, let $\Pc_{m,G}$ denote the set of all in-flow measures $F \in \Pc(\Tc_d \times \Xc)$ that satisfies Assumption \ref{asm: lipschitz of F} and demand constraint given in (\ref{eq: demand profile = marginal dist-x}). Then, by Radon–Nikodym theorem (see e.g. Theorem 32.2 of \cite{billingsley2012probability}) any $F \in \Pc_{m,G}$ admits a probability density function denoted by $f$. Further, let $\Mc_{m,G}$ denotes the set of all positive dis-aggregated in-flow measures $E \in \Pc(\Tc_d \times \Xc \times \Tc_a)$ such that  $F =  \Fc (E) \in \Pc_{m,G}$.
% Further, let $\Mc_{m,G}$ denotes the set of all positive dis-aggregated in-flow measures $E \in \Pc(\Tc_d \times \Xc \times \Tc_a)$ that satisfies the demand constraint given in (\ref{eq:Edependentm}) and such that  $F =  \Fc (E)$ satisfies Assumption \ref{asm: lipschitz of F}, i.e. $F \in \Pc_{m,G}$.

Considering the definition of $\Pc_{m,G}$, the system avoids having a mass of departures at the same time (the in-flow cannot be Dirac). For small values of $G$, an in-flow measure $F \in \Pc_{m,G}$ is very smooth (without any drastic change in a short interval of time). However, when $G$ gets larger, the feasible in-flow measures may have larger fluctuations. Note that Assumption \ref{asm: lipschitz of F} is consistent with the regularity assumption made on the demand profile $m$, i.e., Assumption \ref{asm:regularity}. We have the next proposition about the dynamics of the characteristic travel distance. 
% Furthermore, $\Pc_{m,G}$ is restricted to the in-flow measure that are matched with demand constraint given in (\ref{eq: demand profile = marginal dist-x}).

\begin{proposition}
    \label{prp: congestion dynamics}
    Consider a traffic system with speed function $V$ and in-flow measure  $F \in \Pc_{m,G}$. Then, the characteristic travel distance of the system $z_F$ is the solution of the following set of equations,
    \begin{align}
        \label{eq: system dynamics}
        \begin{cases}
            z_F(t) = \int_0^t V\Big(F\big(S_s(z_F)\big)\Big) ds,
            \\
            S_t(z_F) := \big\{(\tau, \xi)\ \big|\ \tau \in [0,t] \cap \Tc_d,\ \xi \in \big(z_F(t)-z_F(\tau),\infty \big) \cap \Xc\big\}.
        \end{cases}
    \end{align}
\end{proposition}
For the proof refer to appendix \ref{p:prp1}. Note that in the set of equations defined in (\ref{eq: system dynamics}), we use subscript $F$ to emphasize the dependency of the variables on the in-flow measure.
\begin{figure}
\centering
\resizebox*{7.8cm}{!}{\includegraphics{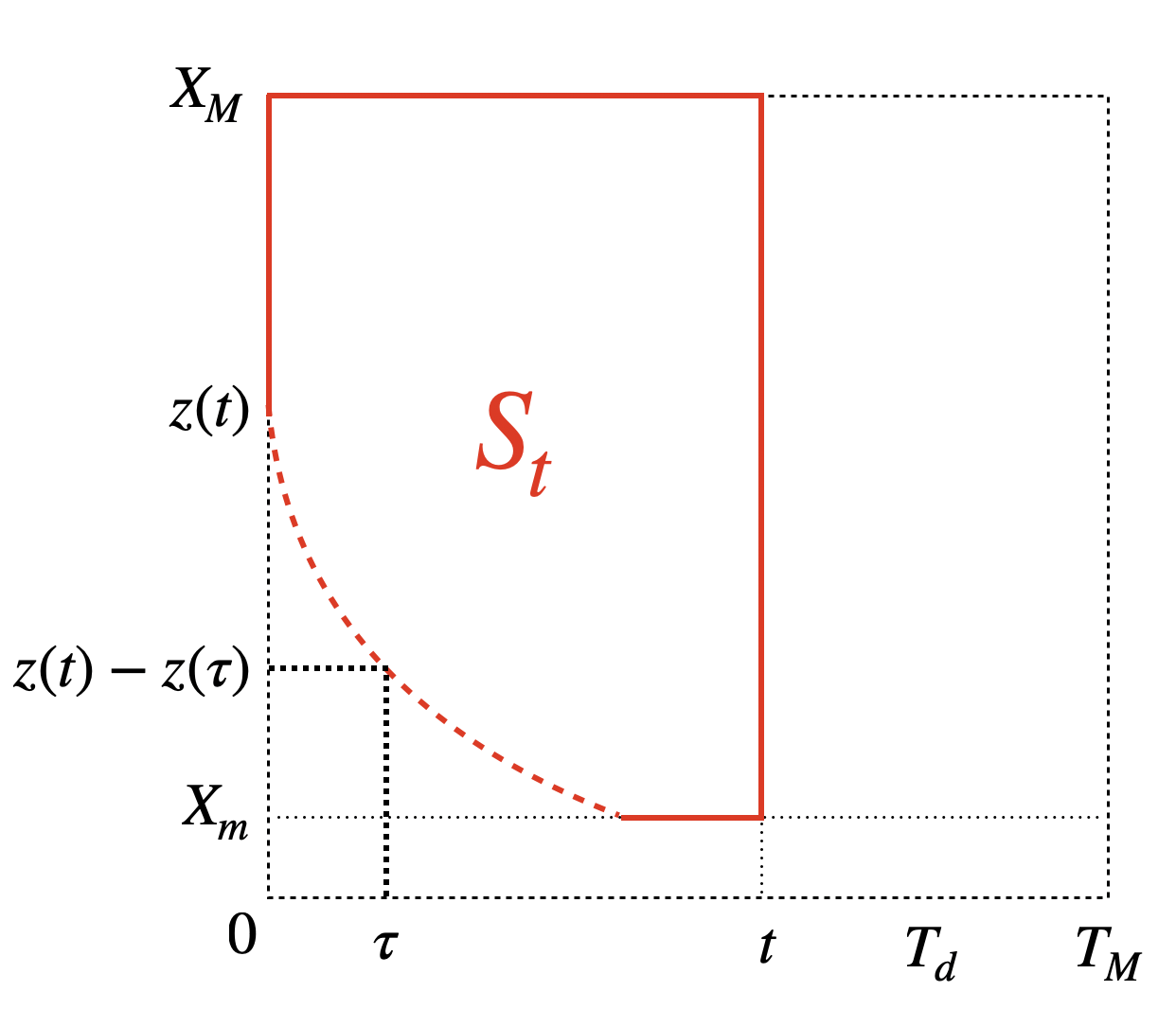}}\hspace{5pt}
\caption{$S_t$ includes the area inside the red lines. The dashed red curve indicates a sample path of $z(t)-z(\tau),\ \tau \in [0,t]$ which is not included in $S_t$.}
\label{fig:S_t}
\end{figure}

Proposition \ref{prp: congestion dynamics} provides the relation between the characteristic travel distance $z_F$ and the in-flow measure $F$, i.e. the distribution of the departures. In fact, $S_t(z_F)$ contains the pairs of the departure times and trip lengths of the users that are travelling at time $t$, in a traffic system with the characteristic travel distance $z_F$. In other words, for all $(t_d,x) \in S_t$, an agent with departure time $t_d$ and trip length $x$ is in the system at time $t$. On the other hand, if $(t_d,x) \notin S_t$, the agent has either not departed or finished her travel before $t$, as illustrated in Figure \ref{fig:S_t} \footnote{The Borel measurability of $S_t$ is obvious.}. 

Note that the set of equations given in (\ref{eq: system dynamics}) should be solved simultaneously. Therefore, we should investigate the existence and uniqueness of the characteristic travel distance derived in Proposition \ref{prp: congestion dynamics}. To address this problem, we need to introduce some notations.
For any compact subset of $\R^n$ such as $C$, $\Cc(C)$ represents the space of all real valued continuous functions defined on $C$. We assume that $\Cc(C)$ is equipped with the uniform norm, i.e.,
\begin{align*}
    \forall u \in \Cc(C):\ \|u\| := \sup_{t \in C} |u(t)|.
\end{align*}
Also, for a constant $M>0$, we define the following norm on $\Cc(C)$ which is  equivalent to the uniform norm on the compact space $C$,
\begin{align*}
    \forall u \in \Cc(C):\ \|u\|_M := \sup_{t \in C} |e^{-tM}u(t)|.
\end{align*}
We denote by $d(\cdot,\cdot)$ and $d_M(\cdot,\cdot)$ the distances associated to $\|\cdot\|$ and $\|\cdot\|_M$, respectively. Also, we define the function $\Uc: \Cc(\Tc) \times \Pc_{m,G} \mapsto \Cc(\Tc)$ such that
\begin{align}
    \label{eq: virtual user fixed point mapping}
    \begin{cases}
    \Uc(z, F) = \Tilde{z},
    \\
    \Tilde{z}(t) = \int_0^t V\Big(F\big(S_s(z)\big)\Big) ds.
    \end{cases}
\end{align}
As demonstrated in the next proposition, systems with smooth in-flow measures, in the sense of Assumption \ref{asm: lipschitz of F}, and Lipschitz continuous speed functions, admit a unique characteristic travel distance.  
\begin{proposition}
\label{prp:existence of z}
For all in-flow measures $F \in \Pc_{m,G}$ and Lipschitz continuous speed functions $V$, there exists a unique function $z_F \in \Cc(\Tc)$ which satisfies the set of the equations given in (\ref{eq: system dynamics}).
\end{proposition}
The proof is given in appendix \ref{p:prp2}

Proposition \ref{prp:existence of z} enables us to obtain $z^*$, the solution of the set of equations given in (\ref{eq: system dynamics}) which is the characteristic travel distance of a system with in-flow measure $F$. The next corollary provides a procedure to obtain $z^*$, based on the successive application of the mapping $\Uc$ defined in (\ref{eq: virtual user fixed point mapping}).
\begin{corollary}
\label{crl: sequence of BFPT}
Fix $F\in \Pc_{m,G}$. Then, starting with an arbitrary element $z_0 \in \Cc(\Tc)$, the sequence $z_l$ defined as
\begin{align}
    \label{eq: sequence of dynamic}
    z_l := \Uc(z_{l-1}, F),\ l\geq 1,
\end{align}
converges to $z^*$ which is the solution of (\ref{eq: system dynamics}).
\end{corollary}
See appendix \ref{p:crl1} for the proof.

% Corollary \ref{crl: sequence of BFPT} provides a procedure to obtain $z^*$, the characteristic travel distance of a system, given the departures distributed wrt $F$.
% Now, fix an arbitrary element $z_0 \in \Cc(\Tc)$. Then, for a given sequence of probability measures in $\Pc_{m,G}$ such as $\{F_m\}_{m \in \N}$, define:
% \begin{align}
%     \label{eq: sequence for F_m}
%     z_{m,n} := \Uc(z_{m,n-1}, F_m),\ n\geq 1.
% \end{align}
% Then, by Corollary \ref{crl: sequence of BFPT}, $z_{m,n} \rightarrow z_m^*$ where $z_m^*$ is the solution of the set of  equations (\ref{eq: system dynamics}) for $F_m$. That means $z_m^*$ is the characteristic travel distance in a system having departures with distribution $F_m$. Similarly, take $z_n$ and $z^*$ as the corresponding sequence to the probability measure $F$.

Consider an arbitrary weakly convergent sequence of probability measures $\{F_k\}_{k \in \N}$ in $\Pc_{m,G}$ such that $F_k \Rightarrow F$. Then, Proposition \ref{lemma: closedness} clarifies that the limit probability measures $F$ lies in $\Pc_{m,G}$, too.
\begin{proposition}
\label{lemma: closedness}
For any $G \in \R^+$, $\Pc_{m,G}$ is a closed subset of $\Pc(\Tc_d \times \Xc)$ in the weak convergence topology.
\end{proposition}
For the proof visit appendix \ref{p:prp3}.

In the following proposition we demonstrate that the characteristic travel distance is continuous with respect to the in-flow measure. Suppose that $z_k^*$ is the solution of the set of  equations given in (\ref{eq: system dynamics}) for $F_k$. That means $z_k^*$ is the characteristic travel distance of a system having departures with distribution $F_k$. Similarly, consider $z^*$ as the corresponding solution to the probability measure $F$. Further, suppose that the probability space $\Pc_{m,G}$ and the set of continuous functions $\Cc(\Tc)$ are equipped, respectively, with the weak and uniform convergence topology. 

\begin{proposition}
\label{prp: contiuity of the dynamics}
Suppose that the speed function $V$ is Lipschitz continuous. Then, the characteristic travel distance is continuous wrt to the in-flow measure. In other words, if $F_k \Rightarrow F$ then $z_k^*\rightarrow z^*$.
\end{proposition}

In the proof of Proposition \ref{prp: contiuity of the dynamics}, see appendix \ref{prf:contiuity of the dynamics}, we provide a convergence bound for the limit of the characteristic travel distances. Thus, the solution of the equations given in (\ref{eq: system dynamics}) is also continuous wrt the dis-aggregated in-flow measure $E$, and the next corollary can be considered as a consequence of the continuity of $\Fc$ (recall that $F=\Fc(E)$). 

\begin{corollary}
\label{eq:crlE}
Suppose that the speed function $V$ is Lipschitz continuous. Then, the characteristic travel distance is continuous wrt to the dis-aggregated in-flow measure $E$. 
\end{corollary}

% Let us consider $F_k = \Fc \left(E_k\right)$, $F = \Fc \left(E\right)$ with $E_k, \: E \in \Pc \left( \Tc_d \times \Xc \times \Tc_a \right)$. If $E_k \Rightarrow E$ as $k\rightarrow \infty$, then $z_k^* \rightarrow z^*$.

The proof is given in appendix \ref{p:crl2}.

Note that the characteristic travel distance of the system depends on $E$ be means of $F$, see (\ref{eq: system dynamics}). 

\subsection{Mean Field Games Departure Time Choice Problem (MFGs-DTCP)}

In this section, using the results derived in the previous sections, we provide a DTCP formulation based on the MFGs approach by assuming that the number of travelers goes to infinity in \ref{eq: DTCP}. Recall that the characteristic travel distance is provided in Proposition \ref{prp: congestion dynamics} and its existence as well as its uniqueness is demonstrated in Proposition \ref{prp:existence of z}. Therefore, considering the objective function of an arbitrary player given in (\ref{eq: MFGs-cost function}) and the relation between the travel time and the characteristic travel distance provided in (\ref{eq: trip time}), we can define the MFGs-DTCP as follows,
\begin{align}
        \tag{MFGs-DTCP}
        \label{eq: MFGs-DTCP}
        \min_{t_d \in \Tc_d} J(t_d; x, t_a; F) &= \alpha T(t_d,x) + \beta\big(t_a - t_d - T(t_d,x)\big)_+ + \gamma\big(t_d + T(t_d,x) - t_a\big)_+
        \\
        \notag
        \textit{s.t.}\\
        &\begin{cases}
        \label{eq:Mfgs-DTCP-CHAR}
        z(t) = \int_0^t V\Big(F\big(S_s(z)\big)\Big) ds,
        \\
        S_t(z) := \big\{(\tau, \xi)\ \big|\ \tau \in [0,t] \cap \Tc_d,\ \xi \in \big(z(t)-z(\tau),\infty \big) \cap \Xc\big\},
        \end{cases}
        \\
        \label{eq:MFGs-DTCP-Triptime}
        &\ \ \ T(t_d,x) = z^{-1}\big(x+z(t_d)\big) - t_d.
\end{align}
Note that in the \ref{eq: DTCP} model defined in (\ref{eq:DTCP-C}), all the three relations should be considered simultaneously, since the choice of an arbitrary player affects the system significantly. However, in the \ref{eq: MFGs-DTCP} model, the set of equations given in (\ref{eq:Mfgs-DTCP-CHAR}) can be investigated independently of the travel time identity provided in (\ref{eq:MFGs-DTCP-Triptime}). This is due to the fact that, as the number of players $n \rightarrow \infty$, the impact of a player on the system vanishes. Moreover, note that the \ref{eq: MFGs-DTCP} model considers the system at a macroscopic level. That is, we do not need to follow the states and decisions of finitely many players\footnote{More precisely, any measure zero subset of agents' indices set is negligible.}. Specifically, if finitely many players change their departure times, the relation \ref{eq:Mfgs-DTCP-CHAR} remains unchanged.

Solving the MFG for the MFGs-DTCP problem is difficult. Indeed it is to be expected that the inflow solution to this problem does not satisfy the regularity Assumption \ref{asm: lipschitz of F} and can exhibit dirac-like concentrations. Therefore we introduce in definition \ref{def: MFE} a relaxed problem with relaxation factor $\varepsilon$. From a physical point of view, $\varepsilon$ can be understood as follows: an agent considers his optimum satisfied if his cost lies within $\varepsilon$ of the minimum cost.

%The following definition clarifies the DTUE as the $\varepsilon$-Mean Field Equilibrium ($\varepsilon$-MFE) of the \ref{eq: MFGs-DTCP}. Note that $\varepsilon$ represents the distance between the global optimum and our solution. We consider $\varepsilon \rightarrow 0$ because if $\varepsilon = 0$, we cannot prove the existence of fixed-point, while the fixed-point is defined and can be solved when $\varepsilon >0$.

\begin{definition}
\label{def: MFE}
Given a constant $\varepsilon \geq 0$, $\Fc(E^*)\in \Pc_{m,G}$ (with $E^*\in \Mc_{m,G}$) is an  $\varepsilon$-Mean Field Equilibrium ($\varepsilon$-MFE) for the \ref{eq: MFGs-DTCP}, if the following relation holds,
\begin{align*}
    E^* \big(\big\{(t_d,x,t_a) \in \Tc_d \times \Xc \times \Tc_a \big| J(t_d;x,t_a; \Fc(E^*) ) \leq J(t_d^0;x,t_a; \Fc(E^*) )+\varepsilon,\ \forall t_d^0 \in \Tc_d \big\}\big) = 1.
    %%  F^* \big(\big\{(t_d,x) \in \Tc_d \times \Xc \big| J(t_d,x; F^*) \leq J(t_d^0,x; F^*)+\varepsilon,\ \forall t_d^0 \in \Tc_d \big\}\big) = 1.
    % \int_{\Tc_d\times \Xc}J(t_d,x; \hat{F}) d\hat{F}(t_d,x)
    % = \inf_{F \in \Pc(\Tc_d \times \Xc)} 
    % \int_{\Tc_d\times \Xc}J(t_d,x; \hat{F}) dF(t_d,x).
\end{align*}
\end{definition}
Note that, Mean Field Equilibrium (MFE) is $\varepsilon$-MFE with $\varepsilon = 0$.
Also, $E^*$ can be expressed as the fixed-point of a map $H: \Mc_{m,G} \mapsto \Mc_{m,G}$ defined as,
\begin{align}
    \label{eq: MFGs-nash-eq}
    \begin{split}
        H(\hat{E}) := \Big\{E \in \Mc_{m,G}\big|\ E \big(\big\{(t_d,x,t_a) \in \Tc_d \times \Xc \times \Tc_a \big| J(t_d;x,t_a; \Fc(\hat{E}) ) \leq  J(t_d^0;x,t_a; \Fc(\hat{E}) ) + \varepsilon,  \\ \ \forall t_d^0 \in \Tc_d \big\}\big) = 1 \Big\}.
    \end{split}
\end{align}
The equivalence of the two definitions provided above holds obviously. In the next section, the existence of an equilibrium for the \ref{eq: MFGs-DTCP} will be examined.

\subsection{Existence of the Equilibrium}

In this section, we show that there exists an equilibrium solution  for the \ref{eq: MFGs-DTCP}. To prove the existence, we first need to examine whether the cost function given in  \ref{eq: MFGs-DTCP} is jointly continuous.
% \begin{lemma}
% \label{lemma: closedness}
% For any $G \in \R^+$, $\Pc_{m,G}$ is a closed subset of $\Pc(\Tc_d \times \Xc)$ in weak convergence topology.
% \end{lemma}
% Proof. Let $\{F_m\}_{m \in \N} \subset \Pc_{m,G}$ such that $F_m \Rightarrow F$. Using the definition of $\Pc_{m,G}$ and Portmanteau theorem, Theorem 2.1. of \cite{billingsley2013convergence}, for all open sets $E \subset \Tc_d \times \Xc$, we have:
% \begin{align*}
%     F(E) \leq \liminf_{m} F_m(E) \leq G\lambda_2(E).
% \end{align*}
% Now, note that Lebesgue measure is outer regular in the sense that a measurable set can be approximated by an open set from outside. That means for all $\varepsilon>0$ and all measurable sets $B\subset \Tc_d \times \Xc$, there exists an open set $E$ such that $B \subset E$ and $\lambda_2(E) < \lambda_2(B) + \varepsilon$. Since $B \subset E$ implies that $F(B) \leq F(E)$, we can write:
% \begin{align*}
%     F(B) \leq F(E) \leq G\lambda_2(E) < G\lambda_2(B) + G\varepsilon.
% \end{align*}
% As $\varepsilon>0$ is arbitrary, it yields:
% \begin{align*}
%     F(B) \leq G\lambda_2(B).
% \end{align*}
% Therefore, $F \in \Pc_{m,G}$ and $\Pc_{m,G}$ is closed under weak convergence.\ep

%Further, Proposition \ref{prp: jointly continuity} provides the last prerequisite to state the existence result by admitting the joint continuity of the cost function.
\begin{proposition}
    \label{prp: jointly continuity}
    Suppose that the speed function is bounded from above and below by $V_{max}$ and $V_{min}$, respectively, such that $V_{max} > V_{min}>0$. Then, the cost function of the \ref{eq: MFGs-DTCP}, $J$, is jointly continuous on $(\Tc_d \times \Xc \times \Tc_a \times \Pc_{m,G})$. Further, the continuity of the cost function on $\Tc_d \times \Xc \times \Tc_a$ is Lipschitz.
\end{proposition}
For the proof refer to \ref{p:prp5}.

%The physical interpretation of the hypothesis $V \geq V_{min}>0$ is the following. The model described in Section \ref{sec:SystemDynamics} is normalized. 
Considering that the velocity $V$ is a function of the congestion $c_t$, the condition $V \geq V_{min}>0$ ensures that the network is not saturated even when $c_t=1$, i.e., the total demand is less than the capacity of the network.

Note that Proposition \ref{prp: jointly continuity} demonstrates  the joint continuity of the cost function only on $(\Tc_d \times \Xc \times \Tc_a \times \Pc_{m,G})$. If we assume that the jointly continuity condition is extendable to $(\Tc_d \times \Xc \times \Tc_a \times \Pc(\Tc_d\times \Xc))$, the problem admits a MFE, by Theorem 4.9 in \cite{lacker2018mean}. Otherwise, we have the next proposition on the existence of the equilibrium.

\begin{proposition}
    \label{prp: existence of the eq}
    For an arbitrary $\varepsilon>0$, there exists a constant $G \in \R^+$ such that \ref{eq: MFGs-DTCP} admits an $\varepsilon$-MFE in the probability space  $\Pc_{m,G}$. This means that there exists an dis-aggregated in-flow measure $E^* \in \Mc_{m,G}$ which is the fixed-point of mapping $H$, given in (\ref{eq: MFGs-nash-eq}).
\end{proposition}
The proof is outlined through appendix \ref{p:prp6}

\textbf{Comment:} The data $m$ (the demand) satisfies the regularity condition expressed by Assumption \ref{asm:regularity}. The regularity constant $M_m$ of the data $m$ and the regularity coefficient $G$ of the $\varepsilon$-MFE should be connected. Indeed the proof of the proposition \ref{prp: existence of the eq} in Appendix \ref{p:prp6}, formula (\ref{eq:BoundOnG}), yields the estimate $G \geq \frac{M_m\, {\rm Lip} (J)}{2 \varepsilon}$. Thus $G$ can be chosen proportional to $M_m$. Further, as $\varepsilon \rightarrow 0$, $G \rightarrow \infty$, which suggests that MFG solutions of the \ref{eq: MFGs-DTCP} problem cannot be found in any $\Mc_{m,G}$ space, and that they could exhibit concentrations of departure times on specific instants (Dirac measures).

\section{MFGs Model for the \ref{eq: MFGs-DTCP}}
\label{Sec:4}

In this section, we aim to characterize the departure time user equilibrium (DTUE) for the Mean Field Games model discussed in the previous section. Recall that Proposition \ref{prp: existence of the eq} guarantees the existence of the departure time equilibrium. 

Consider the optimal behavior of an arbitrary player, assuming that the decisions of the other players are known. Specifically, fix a player with the desired arrival time $t_a$ and the trip length $x$ as well as an in-flow measure $F$ as the proxy for the departure times and trip lengths of the other players. By Proposition \ref{prp:existence of z}, this system has a unique characteristic travel distance $z$. Then, based on  Proposition \ref{prp: congestion dynamics}, we have $v_t := V\big(F(S_t)\big)$, which is the velocity of the system at time $t$. Also, note that (\ref{eq: trip time}) can be written as,
\begin{align}
    \label{eq: trip time by integral}
    \int_{t_d}^{t_d + T(t_d,x)} v_tdt = x.
\end{align}
Then, taking derivative with respect to $t_d$ from both sides of the above equality implies\footnote{Inverse Function Theorem ensures the differentiability of $T$ at any point $t_d$ within $(0,T_{max})$, see for example Theorem 7.4 of \cite{protter2012intermediate}.},
\begin{align*}
    \big(1 + \partial_t T(t_d,x)\big)v_{t_d + T(t_d,x)} - v_{t_d} = 0,
\end{align*}
which yields to,
\begin{align}
    \label{eq: dT/eta}
    \partial_t T(t_d,x) = \frac{v_{t_d}}{v_{t_d + T(t_d,x)}} - 1.
\end{align}
Suppose that $\alpha > \beta$. We apply the first order condition of optimality to determine the equilibrium departure time considering the cost function $J$ defined in (\ref{eq: MFGs-cost function}). In the case that $t_a > t_d + T(t_d,x)$, the third term in the cost function is equal to zero, and we get,
\begin{align*}
    \alpha \partial_t T(t_d,x) - \beta(\partial_t T(t_d,x) + 1) = 0.
\end{align*}
Substituting (\ref{eq: dT/eta}) in the above equation, we get,
\begin{align}
    \label{eq:before t_a}
    \frac{v_{t_d^*}}{v_{t_d^* + T(t_d^*,x)}} = \frac{\alpha}{\alpha - \beta}.
\end{align}
Similarly, if $t_a < t_d + T(t_d,x)$, the second term in the cost function is equal to zero and the following equality can be derived by applying the first order condition,
\begin{align}
    \label{eq:after t_a}
    \frac{v_{t_d^*}}{v_{t_d^* + T(t_d^*,x)}} = \frac{\alpha}{\alpha + \gamma}.
\end{align}
Based on (\ref{eq:before t_a}) and (\ref{eq:after t_a}), it is optimal for an agent who arrives before (after) his desired arrival time to choose the departure time such that the ratio of the system velocity at departure and arrival time be equal to $\frac{\alpha}{\alpha - \beta}$ ($\frac{\alpha}{\alpha + \gamma}$). For on-time agents, based on left and right derivatives of the cost function we can get,
\begin{align}
    \label{eq:on-time t_a}
    \frac{\alpha}{\alpha + \gamma} \leq \frac{v_{t_d^*}}{v_{t_d^* + T(t_d^*,x)}} \leq \frac{\alpha}{\alpha - \beta}.
\end{align}
Summarizing relations provided in (\ref{eq:before t_a}), (\ref{eq:after t_a}), and (\ref{eq:on-time t_a}), we can conclude the following proposition about the optimal choice of an arbitrary agent given the distribution of the others' departures. Note that Proposition \ref{prp: FOC} is in accordance with Proposition 2 in \cite{lamotte2018morning}.
\begin{proposition}
\label{prp: FOC}
The optimal departure time $t_d^*$ of a player having desired arrival time $t_a$ and trip length $x$ with cost function $J$, given in  (\ref{eq: MFGs-cost function}), satisfies the following conditions,
\begin{align}
    \label{eq: FOC1}
    \frac{\alpha}{\alpha + \gamma} \leq \frac{v_{t_d^*}}{v_{t_d^* + T(t_d^*,x)}} \leq \frac{\alpha}{\alpha - \beta}.
\end{align}
Further, for an early and late player we have,
\begin{align}
    \label{eq: FOC2}
    \frac{v_{t_d^*}}{v_{t_d^* + T(t_d^*,x)}} =
    \begin{cases}
    \frac{\alpha}{\alpha - \beta}, \quad t_a > t_d + T(t_d,x),
    \\
    \frac{\alpha}{\alpha + \gamma}, \quad t_a < t_d + T(t_d,x).
    \end{cases}
\end{align}

\end{proposition}
% Note that the cost function defined in (\ref{eq: MFGs-cost function}) is not differentiable at $t_a$. Thus, for the sake of simplicity, we assume that $\alpha = 1$ and $\beta = \gamma$. Additionally, to make the problem tractable and provide insights on how players choose their departure times, we assume that the penalty cost for arriving at $t_d^i + T(t_d^i,x_i)$ instead of the desired arrival time is quadratic. Therefore, we modify the cost function as follows:
% \begin{align}
%     \label{eq: quadratic cost function}
%      J(t_d, x; F) = T(t_d,x) + \frac{\gamma}{2}(t_d + T(t_d,x) - t_a)^2.
% \end{align}
% It is easy to check that all the previous results, particularly the existence of the equilibrium, hold true for the new quadratic cost function, too. Now, we state the following proposition about the optimal choice of a generic player.
% \begin{proposition}
% \label{prp: FOC}
% The optimal departure time $t_d^*$ of a player having desired arrival time $t_a$ and trip length $x$ with cost function $J$, given in  (\ref{eq: MFGs-cost function}), satisfies the following equality,
% \begin{align}
%     \label{prp: FOC}
%     \frac{v_{t_d^*}}{v_{t_d^* + T(t_d^*,x)}} \big(1 +\gamma(t_d^* + T(t_d^*,x) - t_a)\big) = 1.
% \end{align}
% \end{proposition}
% Proof. By checking the first order condition of optimality, we have:
% \begin{align*}
%      \partial_t T(t_d,x) + \gamma\big(1+\partial_t T(t_d,x)\big)(t_d + T(t_d,x) - t_a) = 0.
% \end{align*}
% Then, considering (\ref{eq: dT/eta}), we get (\ref{prp: FOC}).  \ep
Note that the cost function $J$ given in (\ref{eq: MFGs-cost function}) is continuous with respect to departure time on a compact set $\Tc_d$, based on Proposition \ref{prp: jointly continuity}. Therefore, there exists a point at which this function is minimized \footnote{The minimum could be achieved on the boundary of $\Tc_d$. But, the cost function includes a term aiming to minimize the difference between desired and effective arrival time. Thus, we assume that the minimum of $J$ satisfies (\ref{eq: FOC1}) and (\ref{eq: FOC2}).}.
% Moreover, note that $t_d + T(t_d,x) = z^{-1}(x+z(t_d))$ is monotone with respect to $t_d$ while the left side of (\ref{prp: FOC})  is not necessarily monotone in the departure time. Therefore, it is not clear whether the solution of (\ref{prp: FOC}) is unique or not.

Suppose $D: \Tc_a \times \Xc \mapsto \Tc_d$ defines a solution to Proposition \ref{prp: FOC}. That is, $D$ maps the desired arrival time and trip length to the departure time $t_d^*$ which satisfies (\ref{eq: FOC1}) and (\ref{eq: FOC2}), specifically,
\begin{align}
    \label{eq: function of optimal departure time}
    t_d^* = D(t_a, x).
\end{align}
Although Proposition \ref{prp: FOC} characterizes the function $D$, but it is not possible to calculate $D$ explicitly. Proposition \ref{prp: m to e} also provides a more detailed characterization of $D$.
% In the rest of the paper when we use $D(t_a, x)$ as an argument in a function, implicitly it should be understood that it is the first component $t_d^*$ of $D(t_a, x)$ which is the actual argument.

Our next goal is to clarify the relation between the demand profile $m$ and the in-flow measure $F$ in terms of $D$. Consider a population of size $n$, and let $t_a^i$ and $x^i$ denote the $i$-th player's desired arrival time and trip length, respectively. Recall that $\{t_a^i,x^i\}_{i=1}^n$ are i.i.d random variables with the distribution $m$. By (\ref{eq: function of optimal departure time}), $D$ determines the optimal departure times of the players, i.e., $D(t_a^i,x^i) = t_d^i$. Then, based on Glivenko-Cantelli Law of Large Numbers (see e.g. Section 3.2.2 of \cite{cardaliaguet2018short}), almost surely and in $L^1$, $F^n := \frac{1}{n} \sum_{i\in[n]} \delta_{t_d^i, x^i}$ converges weakly to $F$, which is the distribution of $(t_d^i,x^i)$\footnote{We assume $D$ is measurable; thus, $\{t_d^i,x^i\}_{i=1}^n$ are i.i.d RVs, too.}. This result shows that the limit of (\ref{eq: in-flow}), $F$, exists and can be derived based on the demand profile $m$. Using a similar discussion we can show that the limit of (\ref{eq: dis-in-flow}), $E$, exists and represents the the distribution of $(t_d^i,x^i, t_a^i)$. 
% Moreover, $F$ can be specified as the transported probability measure of the demand profile $m$ by the function $D$. That is, for all Borel-measurable subsets of $\Tc_d \times \Xc$ like $B$, we have:
% \begin{align}
%     \label{eq: m to F by D}
%     F(B) = m\big(D^{-1}(B)\big),\quad
%     D^{-1}(B) &= \{\eta \in \Tc_a \times \Xc | D(\eta) \in B\}.
% \end{align}

To clarify the relation between in-flow measure $F$ and demand profile $m$, regarding the function $D$, consider disaggregated in-flow measure $E$. Note that $E(\Delta t_d, \Delta x, \Delta t_a)$ indicates the fraction of the trips having departure time in $\Delta t_d$, trip length in $\Delta x$, and desired arrival time $\Delta t_a$. Then, we can state the following proposition.
\begin{proposition}
    \label{prp: m to e}
    Suppose  $D$, which is defined in (\ref{eq: function of optimal departure time}),  is differentiable with respect to $t_a$ and $\partial_t D(t_a,x) >0,\ \forall t_a \in \Tc_a$
    \footnote{Consider two agents with desired arrival times in $t_a^1<t_a^2$ and the same trip length $x$. Let $t_a^1 = D(t_a^1,x)$ be the departure time of the first player. Then, the virtual user travels a distance of $x$ in the time interval $[t_d^1,t_a^1)$. Therefore, since the velocity is positive, it is rational to assume $D(t_a^1,x) < D(t_a^2,x)$.}. Also, assume that the demand profile $m$ satisfies Assumption \ref{asm:regularity} and dis-aggregated in-flow measure $E$ admits a probability density function $e$.
    Then, we have,
    \begin{align}
        \label{eq: f to m}
        e(D(t_a, x), x, t_a) =
        \frac{m(dx,dt_a)}{\partial_t D(t_a,x)}.
    \end{align}
\end{proposition}
For the proof see appendix \ref{p:prp8}.

Note that (\ref{eq: f to m}) provides the relation between the dis-aggregated in-flow measure and the demand profile that is consistent with the constraint given in (\ref{eq:Edependentm}).

\subsection{MFGs system of equations.}
\label{subsec: MFGs system continuous}

In this section, we discuss the MFGs model for the \ref{eq: MFGs-DTCP} which characterizes the equilibrium of the system. Note that the goal of the MFGs analysis is to examine the equilibrium behavior of travelers (i.e., DTUE) not the individual's optimal departure time. On the other hand, based on Corollary \ref{crl: sequence of BFPT}, a generic player would be able to obtain the characteristic travel distance and determine, using Proposition \ref{prp: FOC}, his strategy given the in-flow measure $F$. Therefore, $F$ is the mean field of the \ref{eq: MFGs-DTCP}, i.e., $F$ captures the required information for a generic agent to describe and analyze the system. 
Denoting the actual arrival time by $\bar{t}_a = t_d + T(t_d,x)$, we can summarize the discussions and results provided in the previous sections to derive the mean field games model:
\begin{empheq}[left={\empheqlbrace}]{align}
    \label{eq:MFGS-optimality}
    &\textstyle\frac{\alpha}{\alpha + \gamma} + \1_{t_a > \bar{t}_a}\big(\frac{\alpha}{\alpha - \beta} - \frac{\alpha}{\alpha + \gamma}\big)
    \leq 
    \frac{v_{t_d}}{v_{\bar{t}_a}}
    \leq
    \frac{\alpha}{\alpha - \beta} + \1_{t_a < \bar{t}_a}\big(\frac{\alpha}{\alpha + \gamma} - \frac{\alpha}{\alpha - \beta}\big)
    \ \mbox{with solution}\ t_d = D(t_a, x),
    \\
    \label{eq:MFGs-distribution}
    &\textstyle e(D(t_a, x), x, t_a) =
        \frac{m(dx,dt_a)}{\partial_t D(t_a,x)},\ f\big(D(t_a, x), x) = \int_{\Tc_a} e(D(t_a, x), x, t_a) dt_a
        ,\ F = \int f(t_d,x) dt_d dx,
    \\
    \label{eq:MFGs-chartravl}
    &\textstyle z(t) = \int_0^t V\Big(F\big(S_s(z)\big)\Big) ds,
    \\
    \label{eq:MFGs-set}
    &\textstyle S_t(z) = \big\{(\tau, \xi)\ \big|\ \tau \in [0,t] \cap \Tc_d,\ \xi \in \big(z(t)-z(\tau),\infty \big) \cap \Xc\big\},
    \\
    \label{eq:MFGs-Triptime}
    &\textstyle T(t_d,x) = z^{-1}\big(x+z(t_d)\big) - t_d,
    \\
    \label{eq:MFGs-constraint}
    &\textstyle F \in \Pc_{m,G}.
\end{empheq}

The  \ref{eq: MFGs-DTCP}  model given in (\ref{eq:MFGS-optimality}-\ref{eq:MFGs-constraint}) can be explained as follows.
Suppose that the decision of the players are captured by the in-flow measure $F$. Using equations (\ref{eq:MFGs-chartravl}) and (\ref{eq:MFGs-set}), and in the light of Corollary \ref{crl: sequence of BFPT}, the associated characteristic travel distance $z$ can be obtained.
Then, a generic agent can employ the relation given in (\ref{eq:MFGs-Triptime}) to determine his travel time. Subsequently, the player is able to obtain the optimal departure time based on (\ref{eq:MFGS-optimality}) along with the function $D$. Finally, the demand profile $m$, which is known, will be transferred according to (\ref{eq:MFGs-distribution}) that specifies the relation between the in-flow measure, dis-aggregated in-flow measure and demand profile. Thus, the optimal distribution of the departure times will be derived as a function $\hat{F}$. Now, based on Definition \ref{def: MFE}, the in-flow measure $F$ would be DTUE of the \ref{eq: MFGs-DTCP} if $\hat{F}$ obtained based on the above procedure is equal to the initial in-flow measure $F$ \footnote{Indeed, relations in (\ref{eq:MFGs-distribution}) consider the dependency of the in-flow measure on the dis-aggregated in-flow measure in Definition \ref{def: MFE}.}.

Note that (\ref{eq:MFGS-optimality}) and (\ref{eq:MFGs-distribution}) are the main components of the model. While the former gives the optimal condition for the decision of a generic player, the latter  captures the distribution of the decisions. The rest are required to make a bridge between relations given in (\ref{eq:MFGS-optimality}) and (\ref{eq:MFGs-distribution}).

% \begin{remark}
% Note that, by discussion after Proposition \ref{prp: FOC}, the uniqueness of the function $D$ does not necessarily hold. Consequently, the transported distribution by $D$, the equilibrium, is not unique either. 
% \end{remark}

\begin{remark}
The \ref{eq: MFGs-DTCP} model provided in (\ref{eq:MFGS-optimality}-\ref{eq:MFGs-constraint}) can be extended  to capture user specific coefficients $\alpha$, $\beta$, and $\gamma$ in the cost function where their distributions are given through demand profile $m$. This paves the way to consider heterogeneous user preferences when solving the DTCP problem.
\end{remark}

\subsection{Discrete approximation of the problem} \label{Sec:5}

In this section, we derive the discrete version of the system of equations given in (\ref{eq:MFGS-optimality}-\ref{eq:MFGs-constraint}) to solve the MFGs model numerically. Let $\Delta t$ and $\Delta x$ denote small intervals in the time and space, respectively, such that $\Delta x \geq V_{max} \Delta t$, where $V_{max}$ indicates the maximum of the network free-flow speed. This means that a trip cannot travel more than $\Delta x$ in a time interval $\Delta t$. We denote the time and space discretization as follows,
\begin{itemize}
    \item The time discretization:
    \begin{equation} \label{eq: time discretization}
        (\tau) = [\tau\Delta t, (\tau+1)\Delta t),
    \end{equation}
    
    \item The space discretization 
    \begin{equation} \label{eq: space discritization}
        (\kappa) = [\kappa\Delta x, (\kappa+1)\Delta x).
    \end{equation}
\end{itemize}
All time intervals $[\tau\Delta t, (\tau+1)\Delta t)$ will be denoted hereafter by $(\tau)$. A similar interpretation holds for $(\kappa)$. Note that these definitions are matched with time horizon $\Tc$ and space set $\Xc$ such that the union of all the defined intervals is equal to the corresponding set, that is $\cup (\tau) = \Tc$ and $\cup (\kappa) = \Xc$. Similarly, let $(\tau_d)$ and $(\tau_a)$, respectively, denote the departure and arrival time intervals where the union of $(\tau_d)$ and $(\tau_a)$ is equal to $\Tc_d$ and $\Tc_a$, respectively.

We define an equivalent discrete version of the demand profile $m$ as follows,
\begin{align}
    \label{eq: discrete demand profile}
    \pi(\tau_a,\kappa) := m\big((\tau_a),(\kappa)\big) = \int_{(\tau_a)} \int_{(\kappa)} m(dt_a, dx).
\end{align}

We assume that the velocity of the system is constant in each time interval $\tau$ and it is captured by $v_{\tau}$. 
Then, the discrete analogous of the optimal condition, given in (\ref{prp: FOC}), can be presented as,
\begin{align}
    \label{prp: FOC-discrete}
    \textstyle\frac{\alpha}{\alpha + \gamma} + \1_{\tau_a > \bar{\tau}_a}\big(\frac{\alpha}{\alpha - \beta} - \frac{\alpha}{\alpha + \gamma}\big)
    \leq 
    \frac{v_{\tau_d}}{v_{\bar{\tau}_a}}
    \leq
    \frac{\alpha}{\alpha - \beta} + \1_{\tau_a < \bar{\tau}_a}\big(\frac{\alpha}{\alpha + \gamma} - \frac{\alpha}{\alpha - \beta}\big),
    % \\
    % \frac{v_{\tau_d}}{v_{\tau_d - T(\tau_d,\kappa)}} \big(1 +\gamma(\tau_d + T(\tau_d,\kappa) - \tau_a)\big) = 1,
\end{align}
where $\bar{\tau}_a := \tau_d + T(\tau_d,\kappa)$ is the actual arrival time interval. Here, with an abuse of notation, $T(\tau_d,\kappa)$ is the travel time of an agent having departure time in $(\tau_d)$ and trip length in $(\kappa)$. Suppose that the function $D$ is a solution of  (\ref{prp: FOC-discrete}). That means, $\tau_d^* = D(\tau_a,\kappa)$ is the optimal departure time interval for a traveler having desired arrival time in $(\tau_a)$ and trip length in $(\kappa)$.
Additionally, let $\mu(\tau_d, \kappa, \tau_a)$ indicate  the fraction of departures in time interval $(\tau_d)$ with trip length in $(\kappa)$ having desired arrival time in $(\tau_a)$. Then, similar to (\ref{eq: f to m}), we can capture the relation between the demand profile $\pi$ and $\mu$ by,
\begin{align*}
    \mu(\tau_d, \kappa, \tau_a) = \frac{\pi(\tau_a,\kappa)\Delta t}{D(\tau_a + 1,\kappa) - D(\tau_a,\kappa)}.
\end{align*}
We also define the discrete characteristic travel distance by,
\begin{align}
    \label{eq: discrete characteristic travel distance}
    \zeta(\theta) := \Delta t \sum_{\tau = 0}^{\theta-1} v_{\tau}.
\end{align}
Then, if we denote by $\Gamma_\theta(\zeta)$ the indices corresponding to the agents that are travelling in the interval $\theta$, we can get,
\begin{align*}
    \Gamma_\theta(\zeta) := \big\{(\tau_d, \kappa)\big| \kappa > \zeta(\theta) - \zeta(\tau)\big\}.
\end{align*}
Moreover, the velocity in a system with the discrete characteristic travel distance $\zeta$ in the time interval $(\tau)$, $v_{\tau}$, would satisfy,
\begin{align*}
    v_{\theta} = V\Big(\sum_{(\tau_d, \kappa) \in \Gamma_\theta(\zeta)}p(\tau_d, \kappa) \Big),
\end{align*}
where $p(\tau_d, \kappa) = \sum_{\tau_a} \mu(\tau_d, \kappa, \tau_a)$ is the fraction of departures in $(\tau_d)$ having trip length in $(\kappa)$ independent of the desired arrival time.
Similarly, we define the the travel time of an agent having departure time in $(\tau_d)$ and trip length in $(\kappa)$, $T(\tau_d,\kappa)$. That is,
\begin{align*}
    T(\tau_d,\kappa) := \zeta^{-1}\big(\kappa + \zeta(\tau_d)\big) - \tau_d.
\end{align*}
Here, $\zeta^{-1}$ shows the inverse of the function $\zeta$, defined in (\ref{eq: discrete characteristic travel distance}).

Therefore, the discrete analogous of the MFGs system defined in  (\ref{eq:MFGS-optimality}-\ref{eq:MFGs-constraint}) can be represented as,
\begin{empheq}[left={\empheqlbrace}]{align}
    \label{eq:MFGs-dis-optimality}
    &\textstyle \textstyle\frac{\alpha}{\alpha + \gamma} + \1_{\tau_a > \bar{\tau}_a}\big(\frac{\alpha}{\alpha - \beta} - \frac{\alpha}{\alpha + \gamma}\big)
    \leq 
    \frac{v_{\tau_d}}{v_{\bar{\tau}_a}}
    \leq
    \frac{\alpha}{\alpha - \beta} + \1_{\tau_a < \bar{\tau}_a}\big(\frac{\alpha}{\alpha + \gamma} - \frac{\alpha}{\alpha - \beta}\big)\ \mbox{with solution}\ \tau_d = D(\tau_a, \kappa),
    \\
    \label{eq:MFGs-dis-distribution}
    &\textstyle \mu\big(D(\tau_a,\kappa,\tau_a), \kappa\big) = \frac{\pi(\tau_a,\kappa)\Delta t}{D(\tau_a + 1,\kappa) - D(\tau_a,\kappa)},\ p(\tau_d, \kappa) = \sum_{\tau_a} \mu(\tau_d, \kappa, \tau_a),
    \\
    \label{eq:MFGs-dis-chartravl}
    &\textstyle \zeta(\theta) := \Delta t \sum_{\tau = 0}^{\theta-1} V\Big(\sum_{(\tau_d, \kappa) \in \Gamma_\tau(\zeta)}p(\tau_d, \kappa) \Big),
    \\
    \label{eq:MFGs-dis-set}
    &\textstyle \Gamma_\theta(\zeta) := \big\{(\tau_d, \kappa)\big| \kappa > \zeta(\theta) - \zeta(\tau)\big\},
    \\
    \label{eq:MFGs-dis-Triptime}
    &T(\tau_d,\kappa) := \zeta^{-1}\big(\kappa + \zeta(\tau_d)\big) - \tau_d,
    \\
    \label{eq:MFGs-dis-constraint}
    &\textstyle \sum_{\tau_d} p(\tau_d,\kappa) = \sum_{\tau_a} \pi(\tau_a,\kappa), \ \forall \kappa.
\end{empheq}
The set of equations given in (\ref{eq:MFGs-dis-optimality}-\ref{eq:MFGs-dis-constraint}) can be explained similar to the ones  provided in (\ref{eq:MFGS-optimality}-\ref{eq:MFGs-constraint}). That means the distribution of the players' decisions, $p$, can be treated as the mean-field of the discrete system.
% Let us elaborate the system of equations \ref{eq:MFGs-dis-optimality}-\ref{eq:MFGs-dis-constraint} by considering the discrete characteristic travel distance, $\zeta$, as the mean-field. Suppose that the discrete characteristic travel distance, $\zeta$, is known for all $\tau$. Using (\ref{eq:MFGs-dis-Triptime}), an arbitrary player determines his travel time. Then, the player is able to find his optimal departure time using (\ref{eq:MFGs-dis-optimality}) and obtain the function $D$. Based on (\ref{eq:MFGs-dis-distribution}), $p$ can be acquired as the distribution of the optimal departure times. By solving (\ref{eq:MFGs-dis-chartravl}) and (\ref{eq:MFGs-dis-set}), one derives the discrete characteristic travel distance of the revised system, $\hat{\zeta}$. The equilibrium, associated to the DTUE but based on the discrete characteristic travel distance, is the fixed-point of this procedure, i.e., $\zeta = \hat{\zeta}$.
Suppose that the decision of the players are given by $p$ for all $\tau_d$ and $\kappa$. Using (\ref{eq:MFGs-dis-chartravl}) and (\ref{eq:MFGs-dis-set}), an arbitrary player can derive the discrete characteristic travel distance $\zeta$, and determine his travel time based on (\ref{eq:MFGs-dis-Triptime}). Then, the player could find his optimal departure time using (\ref{eq:MFGs-dis-optimality}) and obtain the function $D$.
Finally, using (\ref{eq:MFGs-dis-distribution}), $\hat{p}$ can be obtained as the distribution of revised departure times wrt the function $D$, which is derived based on (\ref{eq:MFGs-dis-optimality}). Finally, the DTUE is the fixed-point of this procedure, i.e., $p = \hat{p}$.

\subsection{Equilibrium calculation of the \ref{eq: MFGs-DTCP}}

The equilibrium solution for the DTCP cannot be derived directly from the user optimal control conditions but through an iterative solution method \citep{zhong2011dynamic, ameli2021computational}. In this section, we present an algorithm that can be utilized to numerically solve the discrete MFGs framework.

Recall that the $\varepsilon$-MFE of the DTCP is the fixed-point of (\ref{eq: MFGs-nash-eq}). Therefore, we can apply fixed-point algorithms with a similar optimality conditions to calculate the equilibrium point of the problem. 
In the discrete MFGs framework, determining the equilibrium requires obtaining an approximation for $v_{\bar{\tau}_a}$ given in (\ref{eq:MFGs-dis-optimality}). From the travelers point of view, this approximation enables the players to predict the travel costs required to choose the optimal departure times. The prediction model has to take into account the parameters and evolution of the network, which are captured by the generalized bathtub model. Here, to calculate the equilibrium approximation based on the users decisions, we propose a heuristic algorithm. Our heuristic algorithm is based on the variational inequality theory \citep{noor1988general}. The core idea is to use the delay value ($\bar{t}_a - t_a$) to update the departure times in each iteration. In this case, we also consider the travelers mean speed function with the same desired arrival time as the variable to predict the arrival time of the next simulation with respect to their trip length. Indeed, we reschedule the departure time for each traveler based on equations (\ref{eq:MFGs-dis-optimality}-\ref{eq:MFGs-dis-constraint}). 

The original algorithm is detailed in \cite{friesz2019mathematical}. It is proposed for a continuous dynamic assignment model while we use a discrete version of it based on \cite{ameli2020cross}. Note that in each iteration of the algorithm, a proportion of travelers are selected for rescheduling. 
This proportion is equal to the product of the total demand and a step size. The step size is a coefficient between zero and one that is decreasing during the optimization process \citep{ameli2020simulation}. In this study, the step size is fixed to one over the iteration index. We also add a smart selection process (inspired from \cite{sbayti2007efficient}) to the algorithm in order to speed up the convergence. The process sorts all the trips based on their travel cost (\ref{eq: MFGs-cost function}) and then selects the trips with the higher travel costs for the rescheduling process. Note that in all numerical examples, the length of the time interval in the discrete model is considered as one second. 
%Regarding the  discretization of the time interval length in the generalized bathtub model, it is set to one second. %The algorithm is implemented and applied to departure time choice equilibrium problem of the Lyon North network.

%This proportion is determined by a step size which is decreasing during the optimization process. %In other words, in each iteration, the proportion of the total demand is selected for the rescheduling process. 

%%%%%%%%%%%%%%%%%%%%%%%%%%%%%%%%%%%%%%%%%%%%%%%%%%%%%%%%%%%%%%%%%%%%%%%%%%%%%%%%%%%%%%%%%%
\section{Numerical experiments} \label{Sec:6}

In order to examine the efficiency of our MFGs model, we first compare its performance with one of the recently proposed models in the literature. We then apply the MFGs framework to a large-scale test case in order to evaluate its performance and examine how the optimization procedure to determine the DTUE affects the congestion level of the network's real state. This application is the first one in the literature that addresses the departure time equilibrium on a real large-scale network with a large number of users. 

\subsection{Validation of the MFGs framework} \label{Sec:61}
%?? Check the first use of MFD
\cite{lamotte2018morning} used a quadratic function for the network mean speed function (detailed in \cite{lamotte2018congestion}). They use a trip-based macroscopic fundamental diagram (MFD) model \citep{lamotte2018morning, leclercq2017dynamic} and not the generalized bathtub model. However, both approaches share a common ground and produce similar results in terms of the traffic dynamics. Here, we apply our proposed framework with the exact same demand profile, the same parameters for the cost function, including a smooth approximation of the $\alpha$-$\beta$-$\gamma$ preference, modeled by the marginal utility of the time spent at home $h(t) = \alpha$, the marginal utility of time spent at work $w(t) = \frac{2+\gamma - \beta}{2} + \arctan(4(\bar{t}^i_a - t^i_a))\frac{\gamma+\beta}{\pi}$) and the same parameters for the mean speed function, i.e., the network capacity and free flow speed. The description of all simulation parameters are presented in Section 5.1 of \cite{lamotte2018morning}. %The simulation scenario includes 10 families of 400 trips with $\alpha = 1$, $\beta = 0.4 + \frac{0.2K}{9}$ and $\gamma = 1.5 + \frac{k}{9}$ preferences, where $K={0,1, ..., 9}$. 
Note that trip-lengths are uniformly distributed between 0 and 3. The solution method in \cite{lamotte2018morning} is conducted on a day-to-day basis using a selection method inspired by Method of Successive Average (MSA) and an optimization method based on the grid search (detailed in \cite{lamotte2016morning}). \newref{Table}{tab:Benchmark} compares the optimization results of this model with the one proposed in this paper based on MFGs. The results show that, considering the relative cost, the proposed MFGs method outperforms the previous approach by $5.6\%$. Interestingly, the improvements are much more significant considering the computational effort (the number of iterations) required to converge. Specifically, our numerical results show that the proposed MGFs approach requires $87\%$ less number of iterations to converge. The computation time of each iteration for the proposed method is approximately 3.5 times higher than the reference algorithm because of the MFG approximation. Therefore, the algorithm works $54\%$ faster than the reference algorithm.   
%Both solution methods lead to similar results in terms of optimality, which is measured based on the relative cost. Specifically, the proposed MFGs method  outperforms  the previous approach as the relative cost is improved by $5.6\%$. The improvements are very significant if we consider the computational part as the number of iterations to converge is reduced by nearly $87\%$. %Besides, the MFGs framework approaches the DTUE significantly faster than the grid search with small number of iterations. 
Note that the similarity of the solution quality in terms of the average cost and the total travel time shows that the MFGs framework is consistent with the existing method in the literature for the morning commute problem.

\begin{table}[ht]
\centering
\caption{The quality of the equilibrium approximation.}
\begin{tabular}{l|c|c|c|c}
\multicolumn{1}{c|}{Solution method} & \begin{tabular}[c]{@{}c@{}}Total number\\ of iteration\end{tabular} & \begin{tabular}[c]{@{}c@{}}Convergence indicator\\ {[}relative cost{]}\end{tabular} & \begin{tabular}[c]{@{}c@{}}Average\\ cost\end{tabular} & \begin{tabular}[c]{@{}c@{}}Total travel time\\ {[}sec{]}\end{tabular} \\ \hline
MFGs method                                                                & 259                                                                 & 3.37E-03                                                                            & 12.01                                                  & 26984                                                                 \\
\footnotesize{\cite{lamotte2018morning}}                                                                     & 2000                                                                & 3.57E-03                                                                            & 12.66                                                  & 27362                                                                 \\
Improvement {[}\%{]}                                                              & 87.05                                                               & 5.61                                                                                & 5.13                                                   & 1.38                                                                 
\end{tabular}
\end{table}

\begin{figure}[!ht]
\centering
\subfigure[Time series of accumulation]{%
\resizebox*{8cm}{!}{\includegraphics{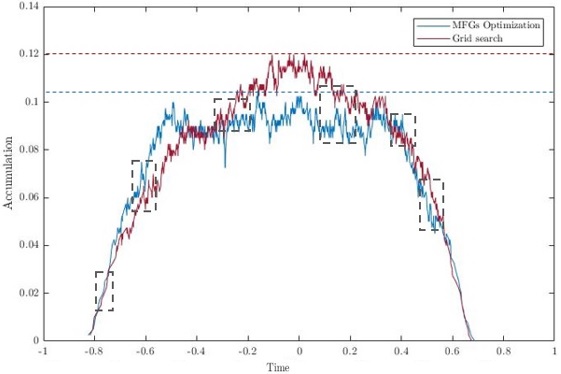}}}
\subfigure[Time series of speed]{%
\resizebox*{8cm}{!}{\includegraphics{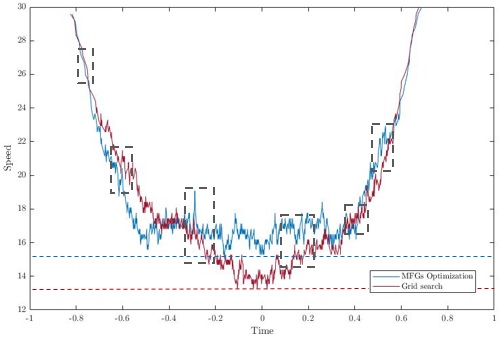}}}
\subfigure[Cumulative departure and arrival curves.]{%
\resizebox*{10cm}{!}{\includegraphics{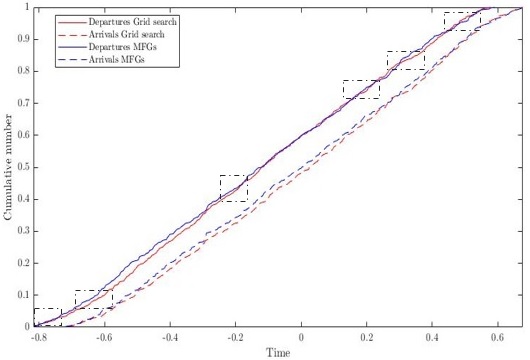}}}
\caption{Simulation results with heterogeneous trip-length: Speed MFD based framework with grid search \citep{lamotte2018morning} versus MFGs framework}
\label{tab:Benchmark}
\end{figure}

To further compare the properties of the final solution for both algorithms, we consider the cumulative departure and arrival curves that provide the characteristics of all trips, and the time-evolution of accumulation and mean speed in the network, see \newref{Figure}{tab:Benchmark}. The MFGs framework provides a solution with a lower maximum accumulation (\newref{Figure}{tab:Benchmark}(a)) and a higher speed (\newref{Figure}{tab:Benchmark}(b)) than the grid search algorithm. It means that the system is closer to the system optimum, defined as the solution where the total travel time of all vehicles is minimum. While it is not the objective function we aim to minimize,  it is interesting to notice that reducing further the total individual costs has a positive impact on the overall system. \newref{Figure}{tab:Benchmark}(c) illustrates  how trips are started  considering trip lengths and departure time. This is a crucial feature as the existing solution methods require prior assumptions on such a sorting to reduce the exploration of the solution space. For example, in \cite{lamotte2018morning}, partial FIFO sorting conditions are mandatory to derive the optimal solutions. Our MFGs framework relaxes such conditions and can provide a full exploration of the solution space. \newref{Figure}{tab:Benchmark} exhibits five time periods (the dash line boxes in each figure) where FIFO patterns are observed in the optimal solution of the grid search solution but no sorting pattern in the MFGs solution. In \newref{Figure}{tab:Benchmark}(c), the inflow rate of the MFGs solution is higher than the grid search while the slope  of the outflow rate is less than the grid search. This test case shows how important it is to relax the sorting assumptions based on the trip lengths to get the optimal solution, which can only be achieved by the proposed MFGs framework.

\subsection{Application of MFGs framework at large scale} \label{Sec:62}

The application of the proposed MFGs framework is easily scalable to much larger instances, which is the main advantage of MFGs over the classic game theory approaches. In this section, we consider a real test case corresponding to the northern part of a metropolis in France (Lyon) and all trips during the morning peak hours, i.e., more than 60,000 trips in total. Note that this section presents the application of our DTUE model to the largest real network with real demand pattern compared to the literature of departure time choice models. %The simulation-based approach allows us to test our model on more complicated configurations. In this section, we scale-up the application of the proposed MFGs model. Here, the goal is to apply our model and assess the optimization process's performance and equilibrium approximation properties in a large-scale network. 

\subsubsection{Description of the test case and demand profile.}

We implement and apply the proposed model to the northern part of Lyon Metropolis (Lyon North). Lyon North network covers 25 $km^2$ and includes 1,883 nodes and 3,383 links. The map is shown in \newref{Figure}{fig:networks}(a). The original demand setting includes all trips during the morning peak hours from 6:30 AM to 10:40 AM (62,450 trips). The data set is published in \cite{HWN8KE_2021} \footnote{\url{https://research-data.ifsttar.fr/dataset.xhtml?persistentId=doi:10.25578/HWN8KE}}. It has been calibrated to represent realistic traffic conditions \citep{krug2019Demand}. All trips have an origin and destination on the real network and a departure times. At the link level of the network (\newref{Figure}{fig:networks}(b)), the origins set contains 94 points and the destinations set includes 227 points. In this study, we only keep the original trip lengths as the generalized bathtub model does not account for the local traffic dynamics. Some trips have origins or destinations outside the covered area (51,215 trips) and will not be considered in the departure time optimization. Note that 11235 trips are fully interior. For those, the original departure time is disregarded and a desired arrival time is assigned. %The network is loaded with travelers of all ODs with a given trip length distribution to represent the morning peak hour (from 6:30 AM to 10:40 AM), based on the study of \cite{krug2019Demand}. The number of trips during this period is 62,450. Some trips start from or end outside the network. Here, we consider only (interior) trips wholly inside the network (11,235 trips) for calculating the departure time equilibrium. 
We divide them into seven classes with different desired arrival times. The desired arrival time of each user is deduced from the real arrival time of the user based on real data \citep{krug2019Demand, alisoltani2019optimal}. The percentage of the trips per class and their desired arrival time are presented in \newref{Table}{tab:Demand}. %The other travelers (51,215 trips) make their trips with given departure time, based on the study of \cite{alisoltani2019optimal}.

\begin{figure}[ht]
\centering
% \subfigure[5 $\times$ 5 grid network]{%
% \resizebox*{6.8cm}{!}{\includegraphics{img/4.jpg}}}\hspace{5pt}
% \subfigure[Manhattan-like network with ring road]{%
% \resizebox*{9.4cm}{!}{\includegraphics{img/5.jpg}}}
\subfigure[Mapping data \copyright Google 2020]{%
\resizebox*{8.5cm}{!}{\includegraphics{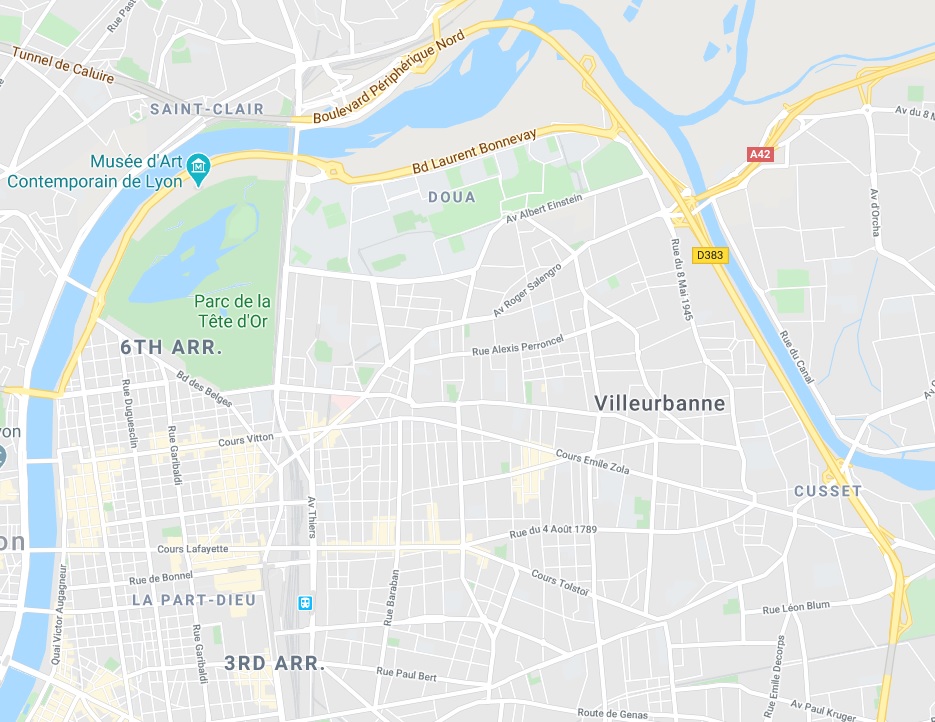}}}\hspace{5pt}
\subfigure[The traffic network using in micro-simulation.]{%
\resizebox*{7.5cm}{!}{\includegraphics{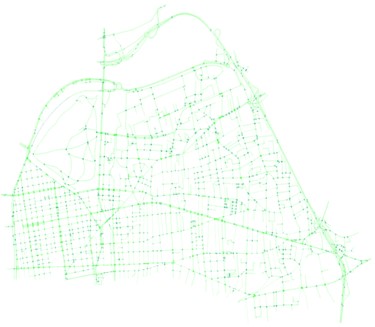}}}
\caption{The northern part of Lyon Metropolis (Lyon North).}
\label{fig:networks}
\end{figure} 

\begin{table}[ht]
\centering
\caption{Demand scenario for Lyon North with multiple desired arrival time}
\begin{tabular}{c|c|c|c|c|c}
\begin{tabular}[c]{@{}c@{}}Class of  \\ the trips\end{tabular} &
\begin{tabular}[c]{@{}c@{}}Share of  \\ the trips\end{tabular} &
\begin{tabular}[c]{@{}c@{}}Number of \\ the trips\end{tabular} &
\begin{tabular}[c]{@{}c@{}}Mean trip \\ length (km)\end{tabular}  &
\begin{tabular}[c]{@{}c@{}}Arrival time \\ interval\end{tabular} & 
\begin{tabular}[c]{@{}c@{}}Estimated Desired \\ arrival time\end{tabular}
\\                                                              \hline
Class 1	    &	13.73\% 	&	1,543 	&	2.53 &  6:30-7:15 &    7:00	\\
Class 2		&	13.84\%    &	1,555	&	2.58 &  7:15-7:45 &    7:30	\\
Class 3	    &	15.42\%   	&	1,732	&	2.55 &  7:45-8:15 &    8:00	\\
Class 4	    &	18.30\%    &	2,056	&	2.65 &  8:15-8:45 &    8:30 \\
Class 5		&	15.05\%    &	1,691	&	2.63 &  8:45-9:15 &    9:00	\\
Class 6	    &	11.82\%   	&	1,328	&	2.70 &  9:15-9:45 &    9:30	\\
Class 7	    &	11.84\%    &	1,330	&	2.63 &  9:45-10:30 &    10:00
\end{tabular}
\label{tab:Demand}
\end{table}

\noindent The network speed function has been calculated in \cite{mariotte2020calibration}. The cost function parameters, i.e., the $\alpha$-$\beta$-$\gamma$ scheduling preferences are defined based on the study of \cite{lamotte2018morning}: $\alpha = 1$, $\beta = 0.4 + \frac{0.2k}{9}$ and $\gamma = 1.5 + \frac{k}{9}$. In order to consider only the heterogeneity of trip length and desired arrival time distributions,  $k$ is fixed to $5$ for all trips in this experiment. 

%\clearpage
%\newpage 

%%%%%%%%%%%%%%%%%%%%%%%%%%%%%%%%%%%%%%%%%%%%%%%%%%%%%%%%%%%%%%%%%%%%%%%%%%%%%%%%%%%%%%%%%%
\subsubsection{Numerical results.} \label{Sec:7}

The optimization process is started with an initial solution where the targeted travelers with a higher trip length in all classes start their trip sooner than others based on the network free-flow speed ($v_{max} = 13.28 m/s$). The heuristic algorithm converges after $56$ iterations to an equilibrium approximation. The results for the convergence pattern is presented in \newref{Figure}{fig:Convergence}(a). The final average cost per traveler is $326.92$, and the figure shows that the final solution is stable. As in the previous test case, the MFGs algorithm converges very fast; however, the performance of heuristic/search algorithms depends on the initial solution. %For instance, \cite{lamotte2018morning} carried out $2000$ iterations for a departure time choice equilibrium with a single desired arrival time.

\newref{Figure}{fig:Convergence}(b) presents the evolution of the network's total travel time during the optimization process. Similar to the convergence pattern, the total travel time decreases and converges to a stable value. Therefore, the final solution can be an equilibrium approximation for the problem. To provide more insights, we also assess the network and equilibrium features overtime during the optimization process.

\begin{figure}[ht]
\centering
\subfigure[Convergence pattern.]{%
\resizebox*{8cm}{!}{\includegraphics{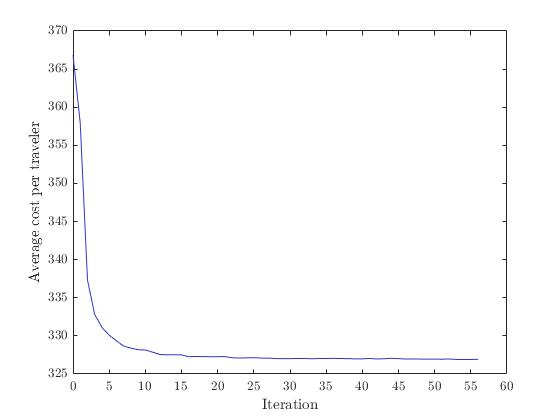}}}\hspace{5pt}
\subfigure[Evolution of the total travel time in the optimization process.]{%
\resizebox*{8cm}{!}{\includegraphics{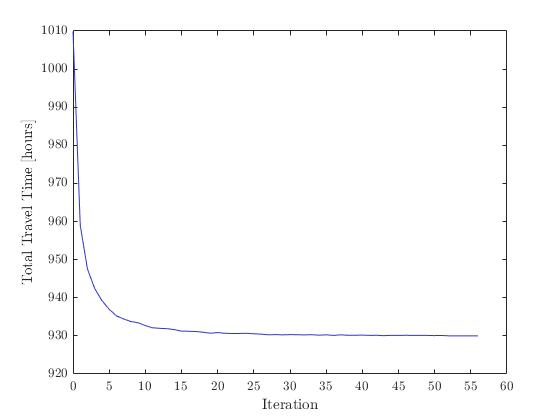}}}
\caption{Results of the optimization process: The travel cost is calculated using (\ref{eq: MFGs-cost function}) and the travel time is the value of $T(t_d, x)$ for each traveler.}
\label{fig:Convergence}
\end{figure}

\begin{figure}[!t]
\centering
\subfigure[Evolution of accumulation for the interior trips: The red line denotes the maximum accumulation of each iteration.]{%
\resizebox*{12cm}{!}{\includegraphics{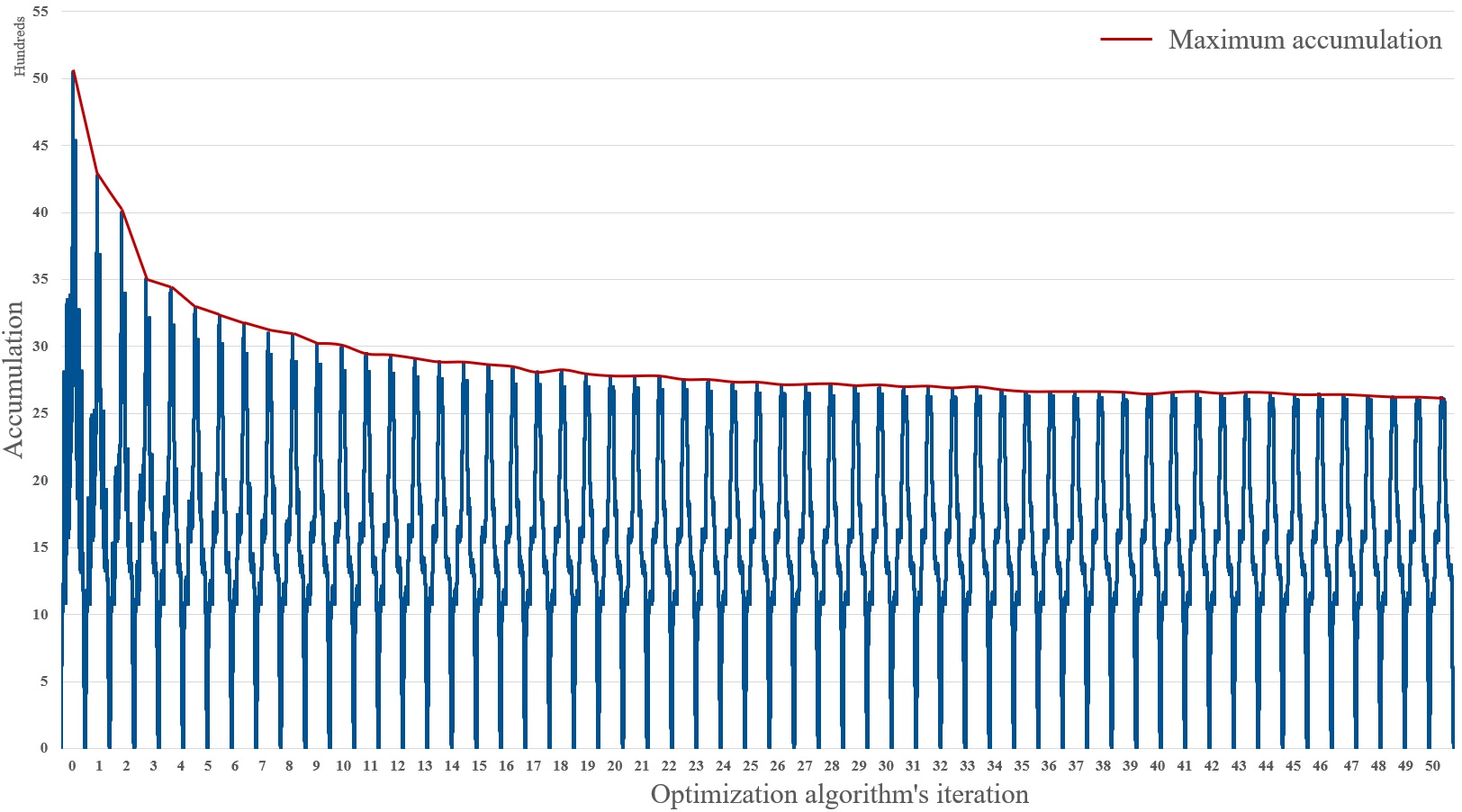}}}\hspace{5pt}
\subfigure[Accumulation of the real state of the network versus equilibrium approximation for all trips.]{%
\resizebox*{12cm}{!}{\includegraphics{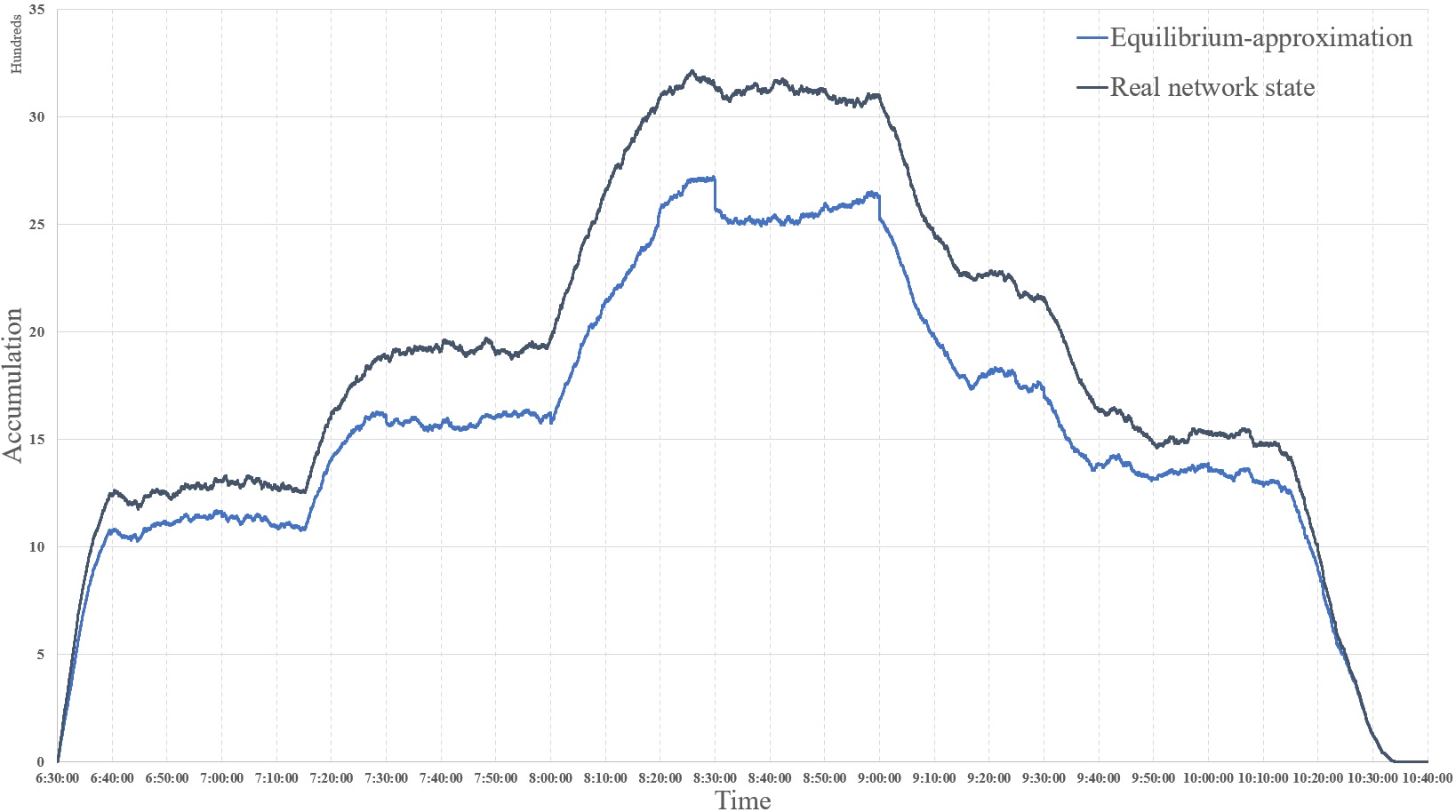}}}
\caption{Results of the network's performance overtime ($\Delta t = 1s$) in the optimization process. Each iteration contains 4.17 hours simulation [6:30 AM - 10:40 AM].}
\label{fig:networkfeatures}
\end{figure} 

\newref{Figure}{fig:networkfeatures} presents the results for the accumulation of the network at each time step ($\Delta t = 1 \: sec$) for the convergence process and the final iteration of the optimization process. In \newref{Figure}{fig:networkfeatures}(a), every blue extrema indicates the evolution of the accumulation in 4.17 hours simulation at one iteration. The curve for the next iteration is started right after the previous one. The results show that the accumulation is also decreasing during the optimization process, and as it is expected, the equilibrium approximation has a low value for the maximum accumulation (red line) of the final solution. Remind that the accumulation evolution in \newref{Figure}{fig:networkfeatures}(a) is drawn for interior trips. % while we have other trips that start or finish their trips outside the network with given departure time. %The shape of the accumulation for each iteration showed that the network reaches its maximal capacity and then unloaded, representing the morning peak of the Lyon North network. 
The equilibrium accumulation for the full demand is shown in \newref{Figure}{fig:networkfeatures}(b). The accumulation time series associated to the original demand patterns with all given departure times is also presented in this figure. This curve is above the solution with optimized departure times. Therefore, the DTUE solution improves the total travel time spent by all users in the system, which is defined by the area between the accumulation time series. %equilibrium approximation, which means the self-optimization of the interior users can improve the network production. Note that these results are obtained while there are other trips that start or finish their trips outside the network with a given departure time.  

The convergence results regarding the different classes of trips are presented in \newref{Figure}{fig:Class}(a). The algorithm's convergence pattern improves in the first three iterations  continuously. However, after the third iteration,  there are small variations for different classes. This is because of our algorithm's heuristic nature that needs to search (explore) the solution space and then exploit it to find a local or the global optimum solution. Note that the exploration rate of heuristic methods depends on the complexity of the solution space and the step size. \newref{Figure}{fig:Class}(b) illustrates another aspect of the equilibrium approximation where each green diamond represents the departure time, and each red circle represents the arrival time of a trip. The duration of a trip is represented by a horizontal blue line between departure and arrival time. The trips of each class are sorted based on their trip length. In \newref{Figure}{fig:Class}(b), the deformations of the distributions for all classes show that non regular sorting pattern matches with the optimal solutions. So, again we show how important it is  to not resort to any prerequisite about the sorting when designing the solution method and defining the optimal conditions. For instance, in \newref{Figure}{fig:Class}(b), the departure and arrival time distributions for class 4 (with desired arrival time 8:30 and the highest demand level) has a deformation on the trip lengths interval [600-850], which illustrates that the partial sorting pattern like FIFO and LIFO does not stand.  %The large-scale results are consistent with the results obtained in the previous section and illustrate that the sorting pattern assumptions do not completely stand for the equilibrium solution.

\begin{figure}[ht]
\centering
\subfigure[Absolute value of the average delay for each class of users.]{%
\resizebox*{8.9cm}{!}{\includegraphics{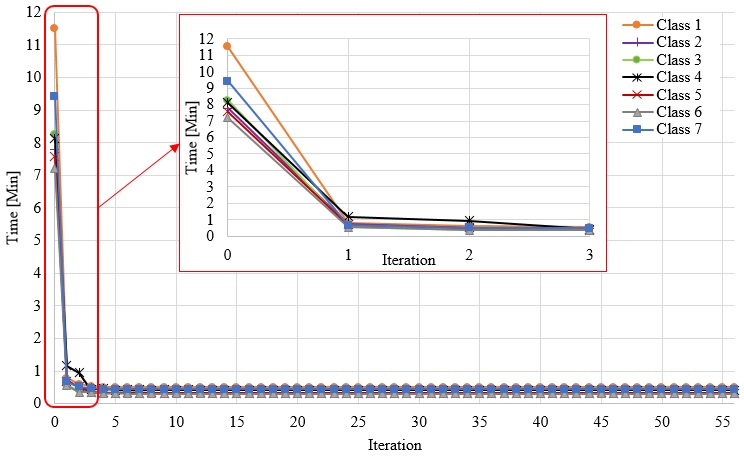}}}\hspace{5pt}
\subfigure[Departure and arrival time distributions of the equilibrium approximation.]{%
\resizebox*{7.1cm}{!}{\includegraphics{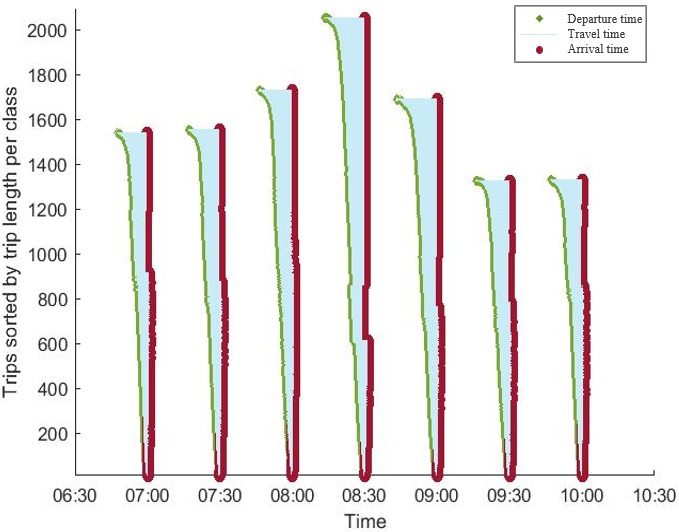}}}
\caption{Optimization results regarding the different classes of trips. Note that there are 11235 interior users in the optimization process.}
\label{fig:Class}
\end{figure} 

%The evolution of the network speed is shown in \newref{Figure}{fig:networkfeatures}(b). Recall that the blue lines around every minima in this figure show the network speed over a full simulation at one iteration. The red line denotes the minimum speed of the system at every iteration. By equilibrating the system, the minimum speed is increasing, which means the network production is increased during the optimization process. 

%Note that the algorithm's convergence pattern for the first five iterations is not improving continuously (e.g., see \newref{Figure}{fig:Convergence}(a)). This is because of our algorithm's heuristic nature that needs to search (explore) the solution space and then exploit it to find a local or the global optimum solution. Remark that the exploration rate of heuristic methods depends on the complexity of the solution space and definition of the step size.

%%%%%%%%%%%%%%%%%%%%%%%%%%%%%%%%%%%%%%%%%%%%%%%%%%%%%%%%%%%%%%%%%%%%%%%%%%%%%%%%%%%%%%%%%%
\section{Conclusion} \label{Sec:Con}
This is the first paper that demonstrates the value of the MFGs approach in departure time equilibrium models. Specifically, this work focuses on modeling and characterizing the departure time choice equilibrium, which is mathematically challenging for large-scale networks, using the  MFGs approach.  We propose a new optimization framework based on the recent findings in transportation systems and game theory. The framework is designed based on mean field game theory coupled with the generalized bathtub model. The MFGs theory allows us to consider a large number of players with different desired arrival times. The idea is that each player in the system optimizes its strategy with respect to the mean-field of the strategies of the other players. Besides, the generalized bathtub model can represent more complex interactions of supply and demand in a transportation system.  

Departure time choices of a group of rational travelers on a traffic network are intrinsically related to how they predict travel time. In this study, we develop a mathematical model through which a generic player predicts the other players' macroscopic behavior. Then, based on this prediction, he derives the dynamics of the system by obtaining the velocity. Having velocity, this player optimizes his departure time strategy. Since the setting is a game with rational players, they look for a Nash equilibrium that can be obtained by a fixed point argument in the procedure of decision making. %To be more precise, the equilibrium is the correct prediction in the sense that it coincides with the resulted distribution by optimizing the departure of a generic player.

Moreover, we  implemented the proposed model for the known setting of the morning commute problem in the literature, and the morning peak hour of the real traffic network of the Lyon North. The numerical results for the first test case demonstrate the value of the MFGs framework compared to existing models in the literature. The large-scale test case shows that the proposed framework is able to represent the equilibrium problem with multiple desired arrival times for a large number of trips that was little addressed before. For the equilibrium calculation, we adapt a heuristic fixed-point algorithm that converges very fast. %compared to the existence of solvers for the departure time choice equilibrium problem. 
The proposed model also provides a good approximation for the equilibrium. The optimization results on both experiments show that optimization based on the mean-field of the users' strategies performs significantly better than the solution methods with Myopia assumption. The gain is much higher when the number of users and the scale of the problem are increased. The equilibrium approximation obtained by the simulation-based optimization contains partial sorting patterns and provide interesting insights on the prevailing sorting assumptions (FIFO or LIFO).  These results are supported by \cite{lamotte2018morning} and underline the importance of the empirical measurements compared to the analytical studies (e.g., \cite{fosgerau2015congestion, daganzo2015distance}). %While the use of trajectory data is increased, it would be nice to investigate the sorting pattern at the microscopic level on the real test case with different demand profiles.     

%For future work, the authors are looking to expand the MFGs framework to consider more complex travel costs (e.g., fuzzy cost) or utility functions (e.g., emissions and tolls) in order to build a more realistic mathematical model. Designing an MFGs framework for dynamic user equilibrium at the microscopic level with a choice model representing users departure time and route choices can be an interesting research direction. Besides, developing optimization algorithms based on the MFGs numerical solution methods looks promising for future works. Regarding the impact of the demand profile on the macroscopic models, the consequences of rapidly varying demand should also be investigated for the DTUE problem.
 
Future researches can be carried out in the following two aspects: (i) analyzing the characteristics of the equilibrium such as uniqueness and stability, which is not easy with heterogeneity of users \citep{lamotte2021monotonicity}; (ii) extend the model by addressing more complex settings and relaxing some of the assumptions. For instance, considering more complex travel costs or utility functions can be an interesting research direction to build a more realistic mathematical model. %The consequences of rapidly varying demand should also be investigated. %Besides, the desired arrival time can be defined as a time window for each class of users, which makes the problem more complex analytically while it can be solved faster. 
Further, developing optimization algorithms based on the MFGs numerical solution methods looks promising for future works.

%%%%%%%%%%%%%%%%%%%%%%%%%%%%%%%%%%%%%%%%%%%%%%%%%%%%%%%%%%%%%%%%%%%%%%%%%%%%%%%%%%%%%%%%%%

% Acknowledgments here
\ACKNOWLEDGMENT{%
This work has received funding from the European Research Council (ERC) under the European Union's Horizon 2020 research and innovation program.
(Grant agreement No 646592 -- MAGnUM project).
We would like to thank Dr. Negin Alisoltani for her helpful discussions and suggestions about the model implementation. 
% Enter the text of acknowledgments here
}% Leave this (end of acknowledgment)

% Appendix here
% Options are (1) APPENDIX (with or without general title) or 
%             (2) APPENDICES (if it has more than one unrelated sections)
% Outcomment the appropriate case if necessary
%
% \begin{APPENDIX}{<Title of the Appendix>}
% \end{APPENDIX}
%
%   or 
%
% \begin{APPENDICES}
% \section{<Title of Section A>}

% \section{<Title of Section B>}
% etc
% \end{APPENDICES}

% References here (outcomment the appropriate case) 

% CASE 1: BiBTeX used to constantly update the references 
%   (while the paper is being written).
\bibliographystyle{informs2014trsc} % outcomment this and next line in Case 1
\bibliography{main.bib} % if more than one, comma separated

% CASE 2: BiBTeX used to generate mypaper.bbl (to be further fine tuned)
% \input{mypaper.bbl} % outcomment this line in Case 2
\begin{APPENDICES}
\section{Proof of Proposition \ref{prp: congestion dynamics}}
\label{p:prp1}
Since $F \in \Pc_{m,G}$, by Radon–Nikodym theorem (see e.g. Theorem 32.2 of \cite{billingsley2012probability}), $F$ admits a probability density function denoted by $f$ Thus, we can employ relation \ref{eq: dynamics of the system} to describe the system congestion. Then, in the light of the equality given in (\ref{eq: virtual user}) and considering that $z'(t) = v_t$, we can rewrite (\ref{eq: dynamics of the system}) as follows,
\begin{equation*}
    \begin{aligned}
    c_t &= \int_0^t \int_0^{X_{max}} f(r,x)dxdr - \int_0^t  \int_0^r v_r f\big(s,z(r)-z(s)\big)ds dr
    \\
    &= \int_0^t \int_0^{X_{max}} f(r,x)dxdr - \int_0^t  \int_s^t v_r f\big(s,z(r)-z(s)\big)dr ds.
    \\
    &= \int_0^t F(dr, \Xc) - \int_0^t F\Big(ds,\big[0,z(t)-z(ds) \big]\Big)
    \\
    &= \int_0^t F\Big(ds,\big(z(t)-z(ds),X_{max} \big]\Big) 
    \\
    &= F\big(S_s(z_F)\big).
\end{aligned}
\end{equation*}
Reusing relation \ref{eq: virtual user} finishes the proof. \ep

\section{Proof of Proposition \ref{prp:existence of z}}
\label{p:prp2}
Consider the mapping $\Uc$ defined in (\ref{eq: virtual user fixed point mapping}). Note that the space $\Cc(\Tc)$ is a Banach space wrt the uniform norm. Therefore, by Banach Fixed-Point Theorem (see Theorem 3.48 of \cite{inbook}), it is sufficient to show that the function $\Uc$ is a contracting mapping, namely there is $\Lambda \in [0,1)$ such that for all $z_1, z_2 \in \Cc(\Tc)$, $d\big(\Uc(z_1, F),\Uc(z_2, F)\big) \leq \Lambda d(z_1,z_2)$.
Hence, consider $z_1, z_2 \in \Cc(\Tc)$ arbitrarily. Using triangle inequality we get:
\begin{align*}
    \big|\Tilde{z}_1(t) - \Tilde{z}_2(t)\big|
    \leq
    \int_0^t \Big|V\Big(F\big(S_s(z_1)\big)\Big) - V\Big(F\big(S_s(z_2)\big)\Big)\Big| dr.
\end{align*}
Then, since $V$ is Lipschitz, we have:
\begin{align*}
    \big|\Tilde{z}_1(t) - \Tilde{z}_2(t)\big|
    &\leq
    \lip(V) \int_0^t \Big|F\big(S_s(z_1)\big) - F\big(S_s(z_2)\big)\Big| ds
    \\
    &\leq
    \lip(V) \int_0^t F\big(S_s(z_1)\ \Delta\ S_s(z_2)\big) ds.
\end{align*}
Using Assumption \ref{asm: lipschitz of F}, we get:
\begin{align*}
    \big|\Tilde{z}_1(t) - \Tilde{z}_2(t)\big|
    \leq
    G\ \lip(V) \int_0^t \int_0^s \big|z_1(s)-z_2(s)\big| +\big|z_1(r)-z_2(r)\big| dr ds.
\end{align*}
Multiplying both sides of the above inequality by $e^{-Mt}$, we obtain:
\begin{align*}
    \big|\Tilde{z}_1(t) - \Tilde{z}_2(t)\big|e^{-tM}
    &\leq
    G\ \lip(V)\int_0^t \int_0^s \big|z_1(s)-z_2(s)\big|e^{-Ms}e^{-M(t-s)} dr ds
    \\
    &+
    G\ \lip(V)\int_0^t \int_0^s 
    \big|z_1(r)-z_2(r)\big|e^{-Mr}e^{-M(t-r)} dr ds
    \\
    &\leq
    2 G\ \lip(V) d_M(z_1,z_2) \frac{Mt - 1 + e^{-Mt}}{M^2},
\end{align*}
where the second inequality is based on the following relations:
\begin{itemize}
    \item $\forall s \in \Tc,\ \big|z_1(s)-z_2(s)\big|e^{-Ms} \leq d_M(z_1,z_2)$,
    
    \item $\int_0^r e^{-M(t-s)} ds = \frac{1-e^{-Mr}}{M}$,
    
    \item $\int_0^t re^{-M(t-r)} dr = \frac{Mt - 1 + e^{-Mt}}{M^2}$,
    
    \item $\int_0^t \frac{1-e^{-Mr}}{M} dr = \frac{Mt - 1 + e^{-Mt}}{M^2}$.
\end{itemize}
Since $\frac{Mt - 1 + e^{-Mt}}{M^2}$ is increasing wrt $t>0$, taking supremum over $t\in \Tc$ yields to:
\begin{align*}
    d_M\big(\Tilde{z}_1,\Tilde{z}_2 \big)
    \leq
    2 G\ \lip(V) \frac{MT_{max} - 1 + e^{-MT_{max}}}{M^2} d_M(z_1,z_2).
\end{align*}
Now, choose $M$ such that $2 G\ \lip(V) \frac{MT_{max} - 1 + e^{-MT_{max}}}{M^2} < 1$. Considering the equivalency between $\|\cdot\|$ and $\|\cdot\|_M$, the proof is complete.\ep

\section{Proof of Corollary \ref{crl: sequence of BFPT}}
\label{p:crl1}
Based on Banach Fixed-Point Theorem (see Theorem 3.48 of \cite{inbook}) and as a result of the proof of Proposition \ref{prp:existence of z}, $z_l$ defined in (\ref{eq: sequence of dynamic}) converges to the fixed point of the function $\Uc$ given in (\ref{eq: virtual user fixed point mapping}).

\section{Proof of Proposition \ref{lemma: closedness}}
\label{p:prp3}
Let $\{F_k\}_{k \in \N} \subset \Pc_{m,G}$ be such that $F_k \Rightarrow F$. Using the definition of $\Pc_{m,G}$ and Portmanteau Theorem (see e.g., Theorem 2.1. of \cite{billingsley2013convergence}), for all open sets $O \subset \Tc_d \times \Xc$, we have:
\begin{align*}
    F(O) \leq \liminf_{k} F_k(O) \leq G\lambda_2(O).
\end{align*}
Now, note that the Lebesgue measure is outer regular in the sense that a measurable set can be approximated by an open set from outside. That means for all $\varepsilon>0$ and all measurable sets $B\subset \Tc_d \times \Xc$, there exists an open set $O$ such that $B \subset O$ and $\lambda_2(O) < \lambda_2(B) + \varepsilon$. Since $B \subset O$ implies that $F(B) \leq F(O)$, we can write:
\begin{align*}
    F(B) \leq F(O) \leq G\lambda_2(O) < G\lambda_2(B) + G\varepsilon.
\end{align*}
As $\varepsilon>0$ is arbitrary, it yields to:
\begin{align}
\label{eq:lipF}
    F(B) \leq G\lambda_2(B).
\end{align}
Additionally, let $Q$ be an arbitrary measurable subset of $\Xc$. For all $k \in \N$, since $F_k \in \Pc_{m,G}$, by (\ref{eq: demand profile = marginal dist-x}), we have,
\begin{align}
    \label{eq:f_k}
    F_k(\Tc_d,Q) = m(Q, \Tc_a), \quad \text{(a.s.)}
\end{align}
where $m$ is the demand profile. On the other hand, by (\ref{eq:lipF}), $F\big(\partial Q \big) \leq G \lambda_2 \big(\partial Q \big) = 0$, and $Q$ is a $F$-continuity set\footnote{$\partial A$, when $A$ is a set, refers to its boundary.}. Then, by the weak convergence of $F_k$ to $F$ and relation (\ref{eq:f_k}), Portmanteau Theorem implies that:
\begin{align*}
    m(Q, \Tc_a) = \lim_k F_k(\Tc_d,Q) = F(\Tc_d,Q),  \quad \text{(a.s.)}.
\end{align*}
Therefore, $F \in \Pc_{m,G}$ and $\Pc_{m,G}$ is closed under weak convergence.\ep

\section{Proof of Proposition \ref{prp: contiuity of the dynamics}}
\label{prf:contiuity of the dynamics}
Let 
$$ \Lambda =  2 G\ \lip(V) \frac{MT_{max} - 1 + e^{-MT_{max}}}{M^2}. $$
Following the proof of Proposition \ref{prp:existence of z}, we know that,
$$ d_M \left( \Uc (z_1,F), \Uc (z_2,F) \right)  \leq \Lambda d_M \left(  z_1,z_2 \right), \quad \forall F\in \Pc_{m,G}, \forall z_1,z_2 \in \Cc (\Tc ),  $$
where $\Uc(z,F)$ denotes the characteristic travel distance of a system having in-flow measure $F$ and primary characteristic travel distance $z$, see (\ref{eq: virtual user fixed point mapping}).
% Now, consider a weakly convergent sequence $F_m \Rightarrow F$ in $\Pc_{m,G}$, and let us denote $z_m^*$, $z^*$ the respective solutions of the equation \ref{eq: system dynamics} for $F_m$, $F$.
Then, we have,
\begin{align*}
   &d_M \left(  z_k^*,z^* \right)  =  d_M \left( \Uc (z^*,F), \Uc (z_k^*,F_k) \right)  \\
   \ \qquad  &\leq   d_M \left( \Uc (z^*,F_k), \Uc (z_k^*,F_k) \right) + d_M \left( \Uc (z^*,F), \Uc (z^*,F_k) \right)  \\
   \ \qquad  &\leq    \Lambda   d_M \left(  z^*,z^*_k \right) + \sup_{0\leq t \leq T_{max}} \exp^{-Mt} \left| \int_0^t \left[ V\left( F \left( S_s (z^*) \right)  \right)- V\left( F_k \left( S_s (z^*) \right)  \right) \right] ds  \right|.  
\end{align*}
Considering that $M$ is chosen such that $\Lambda<1$, we obtain the following bound:
$$ d_M \left(  z_k^*,z^* \right) \leq \frac{\lip(V)}{1-\Lambda}  \sup_{0\leq t \leq T_{max}} \exp^{-Mt} \int_0^t \left|   F \left( S_s (z^*) \right)  - F_k \left( S_s (z^*) \right)    \right| ds $$
which leads to
\begin{equation}\label{eq:FconvergenceBound}
    d_M \left(  z_k^*,z^* \right) \leq \frac{\lip(V)}{1-\Lambda}  \int_0^{T_{max}} \left|   F \left( S_s (z^*) \right)  -  F_k \left( S_s (z^*) \right)    \right| ds.
\end{equation}
The rhs of (\ref{eq:FconvergenceBound}) is bounded by $\frac{2{\rm Lip}(V)  GT_{max} \lambda_2( \Tc \times \Xc) }{1 - \Lambda}$ and converges to 0 as $F_k \Rightarrow F$ by Portmanteau Theorem. Then, using Lebesgue dominated convergence, the rhs of (\ref{eq:FconvergenceBound}) converges to 0 as $k\rightarrow \infty$. Finally, considering equivalency between $\|\cdot\|$ and $\|\cdot\|_M$, the proof of the proposition is complete.\ep

\section{Proof of Corollary \ref{eq:crlE}}
\label{p:crl2}
Due to (\ref{eq:FdependentE}), the bound provided in (\ref{eq:FconvergenceBound}) can then be expressed as 
\begin{equation}\label{eq:EconvergenceBound}
    d_M \left( z_k^*,z^* \right) \leq \frac{\lip(V)}{1-\Lambda}  \int_0^{T_{max}}  \left|   E \left( S_s (z^*) \times \Tc_a \right)   -  E_k \left( S_s (z^*) \times \Tc_a  \right)    \right|ds\, ,
\end{equation}
which shows that if $E_k \Rightarrow E$ as $k\rightarrow \infty$, then $z_k^* \rightarrow z^*$ in the uniform norm. \ep

\section{Proof of Proposition \ref{prp: jointly continuity}}
\label{p:prp5}
The cost function $J$ is defined as follows:
\begin{equation*}
 \begin{aligned}
     J(t_d; x, t_a; F) =
     \begin{cases}
     \alpha T(t_d,x) + \beta\big(t_a - t_d + T(t_d,x)\big),\quad t_a \geq t_d + T(t_d,x)
     \\
     \alpha T(t_d,x) + \gamma\big(t_d + T(t_d,x) - t_a\big),\quad t_d + T(t_d,x) > t_a.
    \end{cases}
 \end{aligned}
\end{equation*}
In order to prove the proposition, it is sufficient to establish the continuity of $T$ as a function of $(t_d,x, t_a)$ and $F$. Considering that $T(t_d,x)=z_F^{-1}\left( x+z_F(t_d) \right)$, it suffices to establish the continuity for $z_F$ and $ z_F^{-1}$. Since $V$ is bounded from above by $V_{max}$, $z_F$ is Lipschitz with $\lip(z_F) \leq V_{max}$. From (\ref{eq:EconvergenceBound}), it follows that for any $t,t'\in \Tc$ and any $F,F' \in \Pc_{m,G}$:
\begin{align*}
 |z_F(t) - z_{F'}(t')|    & \leq  |z_F(t) - z_{F'}(t)| + |z_{F'}(t) - z_{F'}(t')| \\
     & \leq  |z_F(t) - z_{F'}(t)| + V_{max} |t-t'| \\
     & \leq  e^{MT_{max}} \frac{\lip(V)}{1-\Lambda}  \int_0^{T_{max}} ds\, \left|   F \left( S_s (z_F) \right)  -  F' \left( S_s (z_F) \right)    \right| + V_{max} |t-t'|. 
\end{align*}
The continuity of $(t,F) \rightarrow z_F(t) $ is thus established and is obviously Lipschitz wrt $t$. Now, as $V$ is bounded from below by $V_{min}$, it holds that for any $F \in \Pc_{m,G}$:
$ |z_F(t) - z_{F}(t')| \geq V_{min} |t-t'| $. Therefore, $z_F$ is invertible and its inverse is Lipschitz with $\lip(z_F^{-1}) \leq 1/V_{min}$. 
To conclude, consider $F,F' \in \Pc_{m,G}$, $x,x' \in \Xc$, and $t=z_F^{-1}(x)$, $t'=z_{F'}^{-1}(x')$. We have:
\begin{align*}
   z_{F'}(t) - z_{F'}(t') = z_{F'}(t) - z_{F}(t) + x - x'.
\end{align*}
Since
\begin{align*}
    | z_{F'}(t) - z_{F'}(t') | \geq V_{min} |t-t'| = V_{min} | z_{F}^{-1}(x) - z_{F'}^{-1}(x') |
\end{align*}
and 
\begin{align*}
  |z_F(t) - z_{F'}(t)| \leq e^{MT_{max}} \frac{\lip(V)}{1-\Lambda}  \int_0^{T_{max}} ds\, \left|   F \left( S_s (z_F) \right)  -  F' \left( S_s (z_F) \right)    \right|, 
\end{align*} 
we get:
\begin{align*}
   | z_{F}^{-1}(x) - z_{F'}^{-1}(x') | \leq \frac{|x-x'|}{V_{min}} + e^{MT_{max}} \frac{\lip(V)}{(1-\Lambda)\,V_{min}}  \int_0^{T_{max}} ds\, \left|   F \left( S_s (z_F) \right)  -  F' \left( S_s (z_F) \right)    \right|. 
\end{align*}
The continuity of $(x,F) \rightarrow z_{F}^{-1}(x)$ is proved and is Lipschitz wrt $x$. Since $J$ is a piecewise linear function of $t_d$, $t_a$ and travel time function $T$, it follows that $J$ is Lipschitz continuous wrt $t_d$, $x$ and $t_a$. Note that the Lipschitz constant of $J$ depends only on $\lip(V)$, $V_{min}$, $V_{max}$, and $\alpha$, $\beta$, $\gamma$. The continuity coefficient for the dependency of $J$ on $F$ is also dependent on $M$, and $\Lambda$, where $\Lambda$ itself depends on the constant $G$.
\ep

\section{Proof of Proposition \ref{prp: existence of the eq}}
\label{p:prp6}
First, note that $\Pc(\Tc_d\times \Xc \times \Tc_a)$ is a convex compact subset of $\textbf{M}(\Tc_d\times \Xc \times \Tc_a)$, the set of signed measures with bounded variation on $(\Tc_d\times \Xc \times \Tc_a)$, see proof of Theorem 4.10 in \cite{lacker2018mean}. Further, one can show easily that, for all $G \in \R^+$, $\Mc_{m,G}$ is a closed subset of $\Pc(\Tc_d\times \Xc \times \Tc_a)$, with an argument similar to the one given in the the proof of Proposition \ref{lemma: closedness}. Since a closed subset of a compact set is compact, $\Mc_{m,G}$ is a compact subset of $\textbf{M}(\Tc_d\times \Xc \times \Tc_a)$, too. The convexity of $\Mc_{m,G}$ is trivial by its definition.

Now, as $\textbf{M}(\Tc_d\times \Xc \times \Tc_a)$ is a locally convex topological vector space, we can apply the fixed-point theorem of Kakutani (see Theorem 8.6 of \cite{alma991006921539705251}) to prove $H$ admits a fixed point. The convexity and compactness of $H$ is clear. It remains to show that for all $E$, $H(E)$ is non empty and $H$ is usc (upper semi-continuous).

i) Consider $\tilde{E}\in \Mc_{m,G}$. We aim to show that $H(\tilde{E})$ is non-empty. Denote $\tilde{F}=\Fc(\tilde{E})$ and $\tilde{J}(t_d;x,t_a;\tilde{F})=J(t_d;x,t_a;\tilde{F})-\min_{t\in\Tc}J(t;x,t_a;\tilde{F})$. The function $\tilde{J}$ is continuous wrt $F$ and Lipschitz continuous wrt $t_d,x,t_a$. Further, we have $\lip (\tilde{J}) \leq 2 \lip (J)$. Indeed, the function $(x,t_a)\rightarrow \min_{t\in\Tc}J(t;x,t_a;F)$ admits the same Lipschitz constant as $J$. Notice that:
\begin{align*}
   \left|  \tilde{J}(t;x,t_a;\tilde{F}) - \tilde{J}(t';x,t_a;\tilde{F}) \right| \leq \lip (J) | t-t' |, \quad \forall t,t' \in \Tc.
\end{align*}
Denote $U = \{  (t_d,x,t_a) |  \tilde{J}(t_d;x,t_a;\tilde{F}) < \varepsilon \}$ and $U_{x,t_a} = \{  t | \tilde{J}(t;x,t_a;\tilde{F}) < \varepsilon \}$. $U$ is an open set. More precisely, if $(t_d,x,t_a)\in U$, then all $(t'_d,x',t'_a)$ such that 
\begin{align*}
     | t_d-t'_d| + | x-x' | + | t_a-t'_a| < \frac{\varepsilon - \tilde{J}(t_d,x,t_a;\tilde{F})}{2 \lip (J)}
\end{align*}
also belong to $U$.
If $t_d$ is in $U_{x,t_a}$, i.e., $\tilde{J}(t,x,t_a;\tilde{F}) < \varepsilon$, then for all $t'_d$ such that 
\begin{align*}
  | t_d-t'_d | < \frac{\varepsilon - \tilde{J}(t,x,t_a;\tilde{F})}{\lip (J) }.
\end{align*}
Thus, for all $(x,t_a) \in (\Xc, \Tc_a)$, $U_{x,t_a}$ has Lebesgue measure greater than $\frac{2 \varepsilon}{{\rm Lip} (J)}$. It follows that $(x,t_a)\rightarrow \lambda (U_{x,t_a})$ is lsc (lower semi-continuous) and bounded from below by $\frac{2 \varepsilon}{{\rm Lip} (J)}$ where $\lambda$ denotes Lebesgue measure on $\R$. 
Now define the function $\nu$ on $(\Tc_d \times \Xc \times \Tc_a)$ by: 
\begin{align*}
     \nu (t_d,x,t_a) = \frac{1_U (t_d,x,t_a)}{\lambda (U_{x,t_a})}.
\end{align*}
This function is lsc and bounded from above by $\frac{{\rm Lip} (J)}{2 \varepsilon}$.

Finally we define $E$ such that
\begin{align*}
    E(dt_d, dx, dt_a) = \nu(t_d,x,t_a) \lambda(dt_d) m(dx, dt_a).
\end{align*}
By construction, $E$ is positive, has total mass 1 and satisfies the constraint (\ref{eq:Edependentm}). Also, the support of $E$ lies in $U$ by construction; hence, $E$ satisfies $E\in H(\tilde{E})$. It remains to be checked that $E\in \Mc_{m,G}$. It suffices to check that $\Fc(E)$ satisfies Assumption \ref{asm: lipschitz of F}. Consider $B\in \Bc (\Tc_d\times\Xc )$. Using Fubini Theorem (see, e.g., Theorem 18.3 of \cite{billingsley2012probability}) and Assumption \ref{asm:regularity}, regularity condition of $m$, we have, 
\begin{align*}
     \Fc(E) (B)   &  =  E(B\times \Tc_a)    \\
        & = \int_{B\times \Tc_a} dt_d \, dm(x,t_a) \, \nu (t_d,x,t_a)  \\
        & \leq \frac{{\rm Lip} (J)}{2 \varepsilon} \, \int_{B\times \Tc_a} dt_d \, dm(x,t_a)  \\
        & \leq \frac{M_m\, \lip (J)}{2 \varepsilon} \lambda_2(B),
\end{align*}
where $\lambda_2$ is Lebesgue measure on $\R^2$. The above calculation yields to an estimate from below of $G$:
\begin{equation} \label{eq:BoundOnG}
    G \geq \frac{M_m\, {\rm Lip} (J)}{2 \varepsilon}
\end{equation}
It should be noted that the constant ${\rm Lip} (J)$ in (\ref{eq:BoundOnG}) does not depend on $\Fc(E) (B)$ but only on the cost function and on the data $V_{min}$, $V_{max}$ (refer to the proof in Appendix \ref{p:prp5}). Thus choosing $G \geq \frac{M_m\, {\rm Lip} (J)}{2 \varepsilon}$ ensures that $E\in \Mc_{m,G}$. This completes the proof that $H(\tilde{E})$ is non-empty.

ii) We next show that $H$ is usc (upper semi-continuous). Note that $J$ depends on $E$ via $\Fc$. In order to simplify the notations, in this paragraph we will write $J(E)$ for $J(\Fc(E))$. We rewrite the definition of $H$ as follows,
\begin{align*}
   H(E) = \left\{ e\in \Mc_{m,G} \, \Big| \, e\left( \tilde{J}(E) \leq \varepsilon \right) = 1 \right\}.
\end{align*}
In order to show that $H$ is usc it is required to prove that for any open set $W \in \Mc_{m,G}$, the set $H^{-1}W=\{ E | H(E) \subset W\}$ is open, see Page 166 of \cite{alma991006921539705251}. Conversely, denoting $W^c$ and $(H^{-1}W)^c$ the respective complements of $W$ and $H^{-1}W$, it suffices to show that if $W^c$ is closed, $(H^{-1}W)^c$ is closed.

Consider a convergent sequence $\{E_n\}_{n\in\N}$ of elements of $(H^{-1}W)^c$, and let $E$ be the limit of this sequence. We now show that $E\in (H^{-1}W)^c$. For all $n\in\N$, $H(E_n) \not\subset W$ and there exists $e_n\in H(E_n)\cap W^c$. By compactness, we can assume, after extracting a sub-sequence, that the sequence $\{e_n\}_{n\in\N}$ converges weakly towards some $e\in W^c \subset \Mc_{m,G}$. It remains to show that $e\in H(E)$.

Now for any $\eta>0$ there exists $N(\eta)$ such that for $n\geq N(\eta)$, $| \tilde{J}(E_n) - \tilde{J}(E) | < \eta$ (uniformly in $\Cc (\Tc_d\times\Xc\times\Tc_d)$). It implies that
\begin{align*}
\{ \tilde{J}(E_n) \leq \varepsilon  \} \subset \{ \tilde{J}(E) \leq \varepsilon + \eta \},\quad n\geq N(\eta),
\end{align*}
and consequently for all $n\geq N(\eta)$,
\begin{align*}
   e_n \left( \tilde{J}(e) < \varepsilon + \eta \right) = 1.
\end{align*}
Since $e_n \underset{n \rightarrow\infty}{\Longrightarrow} e$, it follows by Portmanteau Theorem that
\begin{align*}
e \left( \tilde{J}(e) \leq \varepsilon + \eta \right) \geq \limsup_{n\rightarrow \infty} e_n \left( \tilde{J}(e) \leq \varepsilon + \eta \right), 
\end{align*}
and thus:
\begin{align*}
    e \left( \tilde{J}(e) \leq \varepsilon + \eta \right) = 1, \quad \forall \eta> 0.
\end{align*}
Finally, $\eta \mapsto 1_{ \{\tilde{J}(e) \leq \varepsilon + \eta\} }$ is monotone decreasing; thus, by Monotone Convergence Theorem, Theorem 4.3.2 of \cite{dudley2018real}, we have:
\begin{align*}
   \lim_{\eta\rightarrow 0} e \left( \tilde{J}(e) \leq \varepsilon + \eta \right) = e \left( \tilde{J}(e) \leq \varepsilon \right) = 1.
\end{align*}
Therefore, we proved that $E\in (H^{-1}W)^c$.
\ep

\section{Proof of Proposition \ref{prp: m to e}}
\label{p:prp8}
Fix $t_a \in \Tc_a$ and $x \in \Xc$. Consider $h$ and $l$ as the small changes in time and space, respectively. The demand with desired arrival time in $\Delta t_a :=(t_a - \frac{h}{2}, t_a + \frac{h}{2})$ and trip length in $\Delta x := (x - \frac{l}{2}, x + \frac{l}{2})$ is equal to $m(\Delta x,\Delta t_a)$ which can be approximated by $m(dx,dt_a)hl$.

On the other hand, since $D$ is increasing wrt the desired arrival time, the departure time of the agents with desired arrival times in $\Delta t_a$ and the trip length $x$ is in the interval $\big(D(t_a - \frac{h}{2}, x),D(t_a + \frac{h}{2},x)\big)$. Then, we approximate the fraction of agents having departure time in $\big(D(t_a - \frac{h}{2}, x),D(t_a + \frac{h}{2},x)\big)$, trip length in $\Delta x$, and desired arrival time in $\Delta t_a$ by 
\begin{align*}
e(D(t_a, x), x, t_a) \big(D(t_a + \frac{h}{2},x) - D(t_a - \frac{h}{2},x)\big) l.
\end{align*}
But, we have:
\begin{align*}
    e(D(t_a, x), x, t_a) \big(D(t_a + \frac{h}{2},x) - D(t_a - \frac{h}{2},x)\big) l \approx m(dx, dt_a)hl.
\end{align*}
Letting $h\rightarrow 0$ yields to the desired result. \ep
\end{APPENDICES}

% \begin{APPENDIX}{Proofs}

% \end{APPENDIX}

\end{document}